\documentclass{article}
\usepackage[hyperindex]{hyperref} 
\usepackage{amsrefs, amsmath,stmaryrd, amssymb, graphicx,bm}
\catcode`\@=11
\@addtoreset{equation}{section}

\newtheorem{Theorem}{Theorem}[section]
\newtheorem{Lemma}{Lemma}[section]

\newtheorem{Definition}{Definition}[section]
\newtheorem{Remark}{Remark}[section]
\newtheorem{Notation}{Notation}[section]

\newcommand{\bProof}{{\bf Proof: }}
\newcommand{\qed}{\hfill $\boxempty $ \bigskip }

\newcommand{\weak}{\rightharpoonup}
\newcommand{\weakstar}{\overset{\star}{\rightharpoonup}}

\newcommand{\Cn}[1]{{\bf [C#1]}}

\newcommand{\Om}{\Omega}
\newcommand{\Set}[1]{\left\{ #1 \right\}}

\newcommand{\inb}{\in_\textrm{b}}
\newcommand{\jump}[1]{\llbracket #1 \rrbracket}
\newcommand{\vc}[1]{\bm{#1}}

\newcommand{\Ah}{\mathcal{A}^h}
\newcommand{\Bh}{\mathcal{B}^h}
\newcommand{\Bk}{\mathcal{B}^h}
\newcommand{\Dk}{\mathcal{D}^h}
\newcommand{\Dh}{\mathcal{D}^h}
\newcommand{\Fh}{\mathcal{F}^h}
\newcommand{\pivh}[1]{\Pi_h^V \left[#1 \right]}
\newcommand{\piqh}[1]{\Pi_h^Q \left[#1 \right]}
\newcommand{\nquad}{\negthickspace\negthickspace
\negthickspace\negthickspace}
\newcommand{\nqquad}{\nquad\nquad}

\date{}
\begin{document}

\title{\bf 1D compressible flow with temperature dependent transport coefficients}
\bigskip
\bigskip

\date{\today}

\author{H. K. Jenssen\thanks{
The research of H.\ K.\ Jenssen was supported in part by the NSF-DMS 0539549 CAREER Award.
}
\thanks{
This paper was written as part of the international research program  
on Nonlinear Partial Differential Equations at the Centre for Advanced  
Study at the Norwegian Academy of Science and Letters in Oslo during  
the academic year 2008-09.  The authors   gratefully acknowledge the
hospitality and support of the Centre for Advanced Study where this research 
 was performed.
}
 \and T. K. Karper\thanks{
A part of this study was conducted while the second author was 
visiting Penn State University. The author expresses his gratitude 
for the hospitality of Penn State. 
}
}

\maketitle

\abstract{
\noindent
We establish existence of global-in-time weak solutions to the one dimensional, 
compressible Navier-Stokes system for a viscous and heat conducting ideal polytropic gas 
(pressure $p=K\theta/\tau$, internal energy $e=c_v \theta$), when the viscosity $\mu$ is constant 
and the heat conductivity $\kappa$ depends on the temperature $\theta$ according to 
$\kappa(\theta) = \bar \kappa \theta^\beta$, with $0\leq\beta<\frac{3}{2}$. This choice of degenerate 
transport coefficients is motivated by the kinetic theory of gasses. 

Approximate solutions are generated by a semi-discrete finite element scheme.
We first formulate sufficient conditions that guarantee convergence to a weak solution. 
The convergence proof relies on weak compactness
and convexity, and it applies to the more general constitutive relations 
$\mu(\theta) = \bar \mu \theta^\alpha$, $\kappa(\theta) = \bar \kappa \theta^\beta$, 
with $\alpha\geq 0$, $0 \leq \beta < 2$ ($\bar \mu,\, \bar \kappa$ constants).
We then verify the sufficient conditions in the case $\alpha=0$ and $0\leq\beta<\frac{3}{2}$.
The data are assumed to be without vacuum, mass concentrations, or vanishing temperatures,
and the same holds for the weak solutions.}

\noindent
{\it 2000 Mathematics Subject Classification: Primary 35Q30; 76N10; Secondary 65M12.}


\section{Introduction}
In one space dimension (1D) the Navier-Stokes system for compressible flow of a heat
conducting and viscous fluid takes the following form in  Lagrangian coordinates:
\begin{align}
	\tau_t &= u_x \qquad \qquad\qquad\qquad\qquad\mbox{(mass conservation)}\label{eq:contequation}\\
	u_t  &= \left(\mu \frac{u_x}{\tau} - p\right)_x \qquad\qquad\qquad\mbox{(momentum balance)} \label{eq:momentumeq} \\
	\mathcal E_t + (up)_x &= \left[\frac{1}{\tau}\Big(\kappa \theta_x+\mu uu_x\Big)\right]_x 
	\qquad\qquad\mbox{(energy conservation)}\label{eq:energy}
\end{align}
Here $x$=Lagrangian 
space variable, $t$=time, and the primary dependent variables are specific volume 
$\tau$, fluid velocity $u$, and temperature $\theta$. The specific total 
energy $\mathcal E=e+\frac{1}{2}u^2$, where $e$ is the specific internal energy. The 
pressure $p$, the internal energy $e$, and the transport coefficients $\mu$ (viscosity) and $\kappa$ 
(heat conductivity), are prescribed through constitutive relations as functions of $\tau$ and $\theta$.
The thermodynamic variables are related through Gibbs' equation $de=\theta dS - pd\tau$, 
where $S$=specific entropy. See \cite{ser} for a derivation of the model.

In this article we establish existence of weak, global-in-time solutions to the 
field equations (\ref{eq:contequation})-(\ref{eq:momentumeq})-(\ref{eq:energy}) for given initial and boundary
data. The flow domain is a finite interval $\Om:=(0,L)$ in Lagrangian coordinates. 
The data we consider are standard (see Theorem \ref{theorem:main} for precise regularity
assumptions): $\tau$, $u$, $\theta$ are prescribed initially on $\Omega$, while $u$ and $\theta$ satisfy 
homogeneous Dirichlet and Neumann conditions, respectively, on $\partial \Omega$;
\begin{equation}\label{eq:bc}
	u = 0, \quad \theta_x = 0, \quad \textrm{ on }\partial \Om \times (0,T).
\end{equation}

Our main interest concerns the transport coefficients $\mu$ and $\kappa$. Their dependence
on $\tau$ and $\theta$ will obviously influence the solutions of the field equations as well as the 
mathematical analysis. Even for 1D flows there is a wide gap between the models furnished by 
physical theories, and the models covered by a satisfactory existence theory. 
We focus on the case of gases for which kinetic theory provides constitutive relations
and we consider only ideal, polytropic gases:
\begin{equation}\label{ideal_polytr}
	p=p(\theta, \tau) = K\, \frac{\theta}{\tau} \,,\qquad e=c_v \theta\,,
\end{equation}
where $K$ (specific gas constant) and $c_v$ (specific heat at constant volume) are 
positive constants. We scale $c_v$ to unity such that $\mathcal E=\theta+\frac{1}{2}u^2$.

According to the first level of approximation in kinetic theory the viscosity 
$\mu$ and heat conductivity $\kappa$ are functions of {\em temperature alone}. 
Furthermore,  the functional dependence is the same for both coefficients. 
See Chapman \& Cowling \cite{cc} or Vincenti \& Kruger
\cite{vk} for a thorough discussion of these issues. If the intermolecular 
potential varies as $r^{-a}$, $r$=intermolecular distance, then $\mu$ 
and $\kappa$ are both proportional to a certain power of the temperature:
\[\mu,\, \kappa\, \propto\,\theta^\frac{a+4}{2a}.\]
For Maxwellian molecules ($a=4$) the dependence is linear, while
for elastic spheres ($a\to+\infty$) the dependence is like
$\sqrt{\theta}$. In any case 
\begin{equation}\label{beta}
	\mu=\bar \mu \theta^b\qquad\mbox{and}\qquad \kappa=\bar \kappa \theta^b 
	\qquad\mbox{for some $b\in (\frac{1}{2},+\infty)$,}
\end{equation}
where $\bar \mu$ and $\bar \kappa$ are constants.
In particular, the transport coefficients tend to zero with $\theta$.
The discrepancy mentioned above is illustrated by the fact that, beyond the regime of 
small and sufficiently smooth data \cite{ko}, there is no global-in-time existence result 
currently available for the Navier-Stokes model 
(\ref{eq:contequation})-(\ref{eq:momentumeq})-(\ref{eq:energy}), 
with constitutive relations \eqref{ideal_polytr} and \eqref{beta}. 

To put our main result Theorem \ref{theorem:main} in perspective let us contrast the last 
statement with what is known about existence of solutions to the compressible Navier-Stokes 
system. (For the purpose of this introduction we concentrate on 1D flows and we consider 
only a small selection of the very extensive literature.) 
The seminal work of Kazhikhov \& Shelukhin \cite{ks} treats the full 
one-dimensional Navier-Stokes system (\ref{eq:contequation})-(\ref{eq:momentumeq})-(\ref{eq:energy}) 
with \eqref{ideal_polytr} and constant transport coefficients. 
Building on earlier work by Nash \cite{na}, Kanel \cite{kan}, and Kazhikhov \cite{kaz}, 
global existence and uniqueness of smooth (i.e.\ $W^{1,2}$) solutions are established in 
\cite{ks} for arbitrarily large and smooth data. 
A key ingredient in the proof is the pointwise a priori estimates on the specific volume  which 
guarantee that no vacuum nor concentration of mass occur. 

Much effort has been invested in generalizing this approach to other cases, and
in particular to models satisfying \eqref{beta}. This has proved to be challenging. 
Temperature dependence of the viscosity $\mu$ has turned out to be especially problematic.
On the other hand, one has been able to incorporate various forms of density dependence
in $\mu$, and also temperature dependence in $\kappa$. Dafermos \cite{d}, Dafermos \& Hsiao \cite{dh}
considered certain classes of solid-like materials in which the viscosity and/or the heat conductivity 
depend on density, and where the heat conductivity may depend on temperature. However, the latter
is assumed to be bounded as well as uniformly bounded away from zero.
Kawohl \cite{k} considered a gas model that incorporates real-gas effects that occur in 
high-temperature regimes. In \cite{k} the viscosity depends only on density (or is constant)
and it is uniformly bounded away from zero, while the thermal conductivity may depend on both 
density and temperature. For example, one of the assumptions in \cite{k} is that there are  
constants $\kappa_0,\, \kappa_1>0$ such that $\kappa(\tau,\theta)$ satisfies
$\kappa_0(1+\theta^q) \leq \kappa(\tau,\theta)\leq \kappa_0(1+\theta^q)$ where $q\geq 2$.  
This type of temperature dependence is motivated by experimental results for gases at very high 
temperatures, see Zel'dovich \& Raizer \cite{zr}. None of these results cover the case 
of a degenerate heat conductivity $\kappa(\theta) = \bar \kappa \theta^b$. 

In the case of isentropic flow a temperature dependence in the viscosity translates into a density 
dependence. For some representative works in this direction see \cite{lxy}, \cite{luoxy}, \cite{mv}, 
\cite{vyz}, \cite{yz}, \cite{yzhu}, and references therein. Finally, for the multi-dimensional 
Navier-Stokes equations there has recently been established various existence results 
where the viscosity and/or heat conductivity depends on $\tau$ or $\theta$, see 
\cite{bd} and \cite{f}. These models do not cover the case of a gas with 
constitutive relations \eqref{ideal_polytr} and \eqref{beta}.

\paragraph{Outline} 
The overall approach in the proof of Theorem \ref{theorem:main} is standard:
apriori pointwise estimates on specific volume  and temperature are coupled 
with higher-order integral estimates to provide sufficient compactness to pass to the limit
in an approximation scheme. However, in both parts of the analysis the temperature 
dependence in the transport coefficients raises some new issues. 

To generate approximate solutions we use a semi-discrete finite element scheme. 
This provides an easy proof of well-posedness of the scheme (Section \ref{scheme}) 
while avoiding some of the cumbersome notation of finite difference schemes. 
In Section \ref{convergence} we formulate certain apriori bounds on the approximations 
and show that these are sufficient for convergence to a weak solution (Theorem \ref{theorem:convergence}). 
The proof employs weak convergence and convexity techniques  \` a la Lions \cite{lions} and Feireisl \cite{f}.
The particular proof we use, working entirely in the Lagrangian frame, seems to be new.
The argument applies to ideal polytropic gases with
$\mu\propto\theta^\alpha$ with $\alpha\geq 0$, and $\kappa\propto\theta^\beta$ with $\beta\in[0,2)$, thus
including the ``standard" case of constant transport coefficients, as well as a full range of 
powers predicted by kinetic theory (see \eqref{beta}). More general constitutive relations could 
presumably be included at the expense of more detailed growth conditions.

In Section \ref{sec:proof} we verify the sufficient apriori bounds in the case where 
$\kappa(\theta) = \bar \kappa \theta^\beta$, with $\beta\in[0,\frac{3}{2})$, and $\mu$ is constant.
Again the result covers constant heat conductivity and a range of power laws suggested 
by kinetic theory. However, we have not been able to establish 
sufficient apriori estimates in the case where also  the viscosity depends on temperature. 
(We point out in the proof 
of Lemma \ref{pws_lemma1} where constant viscosity seems to be essential for the argument.)
To establish the apriori bounds we derive both pointwise estimates and
certain energy estimates. To treat the latter we follow Hoff \cite{hoff} and define energy
functionals that monitor certain weighted $H^1$-norms in the solution. At this point the temperature 
dependence in the heat conductivity requires a careful choice of functionals. On the other hand, the assumption 
of constant viscosity allows us to adopt a standard argument, with only minor changes, to obtain 
the necessary pointwise estimates on $\tau$ and $\theta$.
For completeness this part is included in Section \ref{sec:pointwise}.

\section{Main result}
It is convenient to formulate the approximation scheme for the temperature field instead of 
the total energy $\mathcal E$. Consequently we consider a weak form of the temperature equation:
\begin{equation}\label{eq:tempeq}
	\theta_t - \left(\frac{\kappa\theta_x}{\tau}\right)_x = \mu\frac{u_x^2}{\tau} - pu_x\,. 
\end{equation}
In the following definition $\Omega=(0,L)$, with $0<L<\infty$, and we 
assume that $\mu$, $\kappa$, and $p$ are given, smooth, non-negative functions of $(\theta,\tau)$. 
\begin{Definition}[Weak solution]\label{def:weak}
We say that $(\tau, u, \theta)$ is a weak solution of the compressible Navier-Stokes system 
\eqref{eq:contequation}, \eqref{eq:momentumeq}, \eqref{eq:tempeq} on $\Om \times [0,T)$,
with boundary conditions \eqref{eq:bc} and initial data $(\tau_0, u_0, \theta_0)$ 
satisfying
\begin{equation}\label{weak:initial}
	 (\tau_0, u_0, \theta_0) \in L^2(\Om)\times L^2(\Om) \times L^2(\Om)\,,
\end{equation}
provided that 
$\tau \in W^{1,2}(0,T;L^2(\Om))$, $u \in L^\infty(0,T;L^2(\Om))\cap L^2(0,T;W^{1,2}_0(\Om))$, 
$\theta \in L^\infty(0,T;L^2(\Om))\cap L^2(0,T;W^{1,2}(\Om))$, 
\begin{equation}\label{weak:masseq}
	\tau_t = u_x, \quad \text{a.e in $\Om \times (0,T)$}, \quad \tau(x,0) = \tau_0(x),~ \text{ a.e $x \in \Om$},
\end{equation}
and that 
for all $\phi \in C_c^\infty(\Om\times [0,T))$:
\begin{align}
	&\nquad \int_0^T \int_\Om u \phi_t - \Big[\frac{\mu u_x}{\tau} - p\Big] \phi_x\,dxdt 
	+\int_\Om u_0 \phi_0\,dx = 0, \label{weak:momentumeq}\\
	&\nquad \int_0^T \int_\Om \theta \phi_t - \frac{\kappa}{\tau}\theta_x \phi_x + 
	\Big[\frac{\mu u_x}{\tau} - p\Big] u_x\phi\,dxdt 
	+\int_\Om \theta_0 \phi_0\,dx = 0.
	\label{weak:tempeq}
\end{align}
In addition all terms in \eqref{weak:momentumeq} and \eqref{weak:tempeq} 
are required to be integrable. 
\end{Definition}
In the remainder of the paper we restrict ourselves to the case of an ideal polytropic gas, i.e.\ 
\eqref{ideal_polytr} holds. Furthermore, for concreteness, and motivated by kinetic theory, 
we only consider the case where $\mu$ and $\kappa$ are proportional to (possibly different) 
powers of $\theta$:
\begin{equation}\label{mu_kappa}
	\mu(\theta)=\bar \mu \theta^\alpha, \qquad \kappa(\theta)=\bar \kappa \theta^\beta\,,
\end{equation}
where $\bar \mu, \, \bar \kappa > 0$, and $\alpha,\, \beta\geq 0$ are constants.
To state the main result we define the following functionals:
\begin{align*}
	\mathcal{A}(t) &= \sup_{s \in (0,t)}\int_\Om  |u_x|^2~dx
	+\int_0^t \int_\Om  |u_t|^2~dxdt, \\
	\mathcal{B} (t) &= \sup_{s \in (0,t)}\int_\Om L(\theta)~dx 
	+ \int_0^t\int_\Om  |L'(\theta^h)_x|^2~dxdt, \\
	\mathcal{D} (t) &= 
	\sup_{s \in (0,t)} \frac{1}{2}\int_\Om \frac{1}{\tau}|L'(\theta)_x|^2~dx
	+\int_0^t\int_\Om \kappa(\theta)|\theta_t|^2~dxdt\,,
\end{align*}
where 
\begin{equation}\label{L}
	L(z):=\int_0^z\int_0^\xi \kappa(\eta)\, d\eta d\xi\,.
\end{equation}
Our main result is:
\begin{Theorem}\label{theorem:main}
Consider the one-dimensional, compressible Navier-Stokes system \eqref{eq:contequation}, 
\eqref{eq:momentumeq}, and \eqref{eq:tempeq} for an ideal, polytropic gas \eqref{ideal_polytr}.
Assume that the transport coefficients $\mu$ and $\kappa$ satisfy \eqref{mu_kappa} with 
$\alpha=0$ and $0 \leq \beta < \textstyle\frac{3}{2}$. Let the initial data satisfy
\begin{equation}\label{data_strong}
	(\tau_0, u_0, \theta_0)\in L^2(\Om) \times W^{1,2}(\Om) \times W^{1,2}(\Om)\,,
\end{equation}
and be such that
\begin{equation}\label{data_bnds}
C_0^{-1} \leq \tau_0(x) \leq C_0\,,
\qquad C_0^{-1}  \leq \theta_0(x),
\qquad\mbox{for a.a.\, $x \in \Om$,}
\end{equation}
where $C_0>0$ is a constant. 

Then, for any finite time $T>0$ there exists a weak solution $(\tau, u, \theta)$ on 
$\Om \times [0,T)$ of \eqref{eq:contequation}, \eqref{eq:momentumeq}, and 
\eqref{eq:tempeq}, with the boundary conditions \eqref{eq:bc} 
and initial data $(\tau_0, u_0, \theta_0)$.  
Furthermore, there exists $C>0$, depending on the parameters $K$, 
$c_v$, $\bar \mu$, $\bar \kappa$, $\beta$, the initial data, and $T$, such that 
\begin{equation}\label{main:pointwise}
C^{-1} \leq \tau(x,t) \leq C\,, \qquad 	\theta(x,t) \geq C^{-1}\,,\quad\textrm{ a.e in }\Om \times (0,T),
\end{equation}
and
\begin{equation}\label{main:gradient}
	\sup_{t \in (0,T)} \big[\mathcal{A}(t) + \mathcal{B}(t) + \mathcal{D}(t)\big] \leq C\,.
\end{equation}
The weak solution can be obtained as the pointwise a.e.\ limit in $\Om \times [0,T)$ of
solutions to the semi-discrete finite element scheme described in Definition \ref{def:fem}.
\end{Theorem}
The proof is detailed in the following sections and summarized in Section \ref{biggie}.
\begin{Remark}
	Concerning the regularity of the initial data we note that the result in \cite{hoff} 
	covers discontinuous data  under the assumption of BV regularity of $\tau_0$ and $u_0$. 
	The scheme we consider (Definition \ref{def:fem}) is well-defined for 
	$(\tau_0, u_0, \theta_0)\in \big[L^2(\Om)\big]^3$; the higher regularity in \eqref{data_strong}
	is required to bound the initial values of $\mathcal{A}$, $\mathcal{B}$, $\mathcal{D}$.
	We also mention that \eqref{main:pointwise} may be refined to 
	give {\em time-independent} density bounds by adopting the techniques in \cite{cht}.
\end{Remark}

\section{Finite element scheme}\label{scheme}
We define a semi-discrete finite element scheme 
approximating the compressible Navier-Stokes system.
We verify the basic bounds for mass, energy and entropy for the scheme,
and show that the scheme is well-defined globally in time.

\subsection{Approximation scheme}
Let $\{\mathbb{E}_h\}_{h>0}$ be a family of uniform meshes of $\Om$, where 
$h$ is the mesh size and it is assumed that $N= \frac{|\Omega|}{h}\in\mathbb N$.
For each such $h$, the vertices of $\mathbb{E}_h$ are $x_i = ih$, $i=0, \ldots, N$.
For each interval $E\in \mathbb{E}_h$, $\mathbb{P}_k(E)$ denotes the space of 
polynomials of maximal order $k$ on $E$.
On $\mathbb{E}_h$ we define the space of piecewise constants,
$$
	Q_h(\Om) = \Set{\phi \in L^2(\Om); \phi|_E \in \mathbb{P}_0(E), \ \forall E \in\mathbb{E}_h}\,,
$$
and the space of continuous piecewise linears,
$$
	V_h(\Om) = \Set{v \in C(\Om); v|_E \in \mathbb{P}_1(E), ~\forall E \in\mathbb{E}_h}\,.
$$
Elements in $Q_h(\Om)$ are taken to be continuous from the right, as are derivatives of 
elements in $V_h(\Om)$. In order to incorporate Dirichlet boundary conditions we define the space 
$$
	V_h^0(\Om) = \Set{v \in V_h(\Om); v|_{\partial \Om} = 0}\,.
$$
For $Q_h(\Om)$ we define the projection operator 
$\Pi_h^Q: L^p(\Om) \rightarrow Q_h(\Om)$, $p\in[1,\infty]$ by
\begin{equation}\label{proj1}
	\int_E [\Pi_h^Q \phi](x) ~dx = \int_E \phi(x) ~dx, \quad \forall E \in\mathbb{E}_h\, .
\end{equation}
For $V_h(\Om)$ we define the projection operator 
$\Pi_h^V: W^{1,2}(\Om) \to V_h(\Om)$ by
$$
	[\Pi_h^V \phi](x_i) = \phi(x_i), \quad i=0, \ldots, N\,,
$$
such that 
\begin{equation}\label{proj2}
	\big(\Pi_h^V \phi\big)_x = \Pi_h^Q (\phi_x)\,.
\end{equation}
For any $q\in Q_h(\Om)$ we set $\jump{q}_i := q|_{E_{i+1}} - q|_{E_i}$, where $E_i$ denotes the $i$th 
element of $\mathbb{E}_h$.
\begin{figure}
	\begin{center}
	\includegraphics[width= 0.5\textwidth]{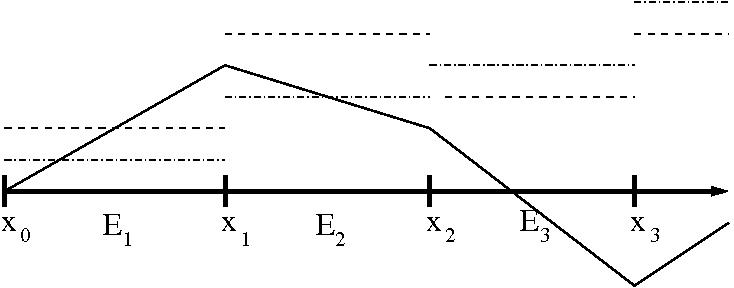}
\end{center}
\caption{A solution of the scheme at a fixed time. The specific volume $\tau^h$
and temperature $\theta^h$ are piecewise constant while the 
velocity $u^h$ is piecewise linear. 
Note that the nodes are numbered from $0$ while the elements are numbered from $1$.}
\label{figure:grid}
\end{figure}
\begin{Definition}[Semi-discrete finite element scheme]\label{def:fem}
Let the constitutive relations \eqref{ideal_polytr} and \eqref{mu_kappa} hold and fix $h>0$ with 
$N:= \frac{|\Omega|}{h}\in\mathbb N$.
For initial data $(\tau_0, u_0, \theta_0)$ satisfying \eqref{weak:initial} we set 
\begin{equation}\label{data_h}
	(\tau^h_0, u^h_0, \theta^h_0) := (\piqh{\tau_0}, \Pi_h[u_0], \Pi_h^Q[\theta_0]),
\end{equation}
where $\Pi_h: L^p(\Om) \rightarrow V^h(\Om)$ is the $L^2$-projection.
For a given time $T>0$ we determine functions 
\begin{equation*}
	(\tau^h, u^h, \theta^h)(t) \in Q_h(\Om)\times V_h^0(\Om)\times V_h(\Om), \qquad t \in (0,T),
\end{equation*}
such that
\begin{align}
	&\int_\Om \tau^h_t \phi^h~dx = \int_\Om u^h_x \phi^h~dx 
	\qquad\quad\forall \phi^h \in Q_h(\Om) \label{fem:contequation} \\
	&\int_\Om u^h_t v^h~dx = -\int_\Om \Big[\frac{\mu(\theta^h)u^h_x}{\tau^h} - p(\theta^h, \tau^h)\Big]v^h_x~dx
	\qquad\quad\forall v^h \in V_h^0(\Om) \label{fem:momentumeq} \\
	&\int_\Om \theta^h_t\psi^h~dx 
	+ \frac{1}{h}\sum_{i=1}^{N-1}G_i(\tau^h)\jump{L'(\theta^h)}_i\jump{\psi^h}_i \nonumber \\
	&\qquad\qquad= \int_\Om \Big[\frac{\mu(\theta^h)|u^h_x|^2}{\tau^h} - p(\theta^h, \tau^h)u^h_x\Big]\psi^h ~dx
	 \qquad\quad\forall \psi^h \in Q_h(\Om),\label{fem:tempeq}
\end{align}
where $G_i(\tau^h)= 2\big(\tau^h|_{E_i} + \tau^h|_{E_{i+1}}\big)^{-1}$ and 
$L(z)$ is given by \eqref{L}.
\end{Definition}

\subsection{Basic estimates and well-posedness of the scheme}
To show that the scheme is well-defined we recast \eqref{fem:contequation}-\eqref{fem:tempeq} 
as an ODE on the (finite dimensional) finite element space $Q_h(\Om) \times V_h^0(\Om) \times Q_h(\Om)$. 
Assuming for now  that a solution $(\tau^h, u^h, \theta^h)$ exists, we  show  
mass and energy conservation, and entropy balance, for the scheme. 
With $\phi^h = 1$ equation \eqref{fem:contequation} gives
\begin{equation}\label{eq:massconv}
	\int_\Om \tau^h(x,t) ~dx \equiv \int_\Om \tau^h_0(x) ~dx\,.
\end{equation}
Next, \eqref{fem:momentumeq} with $v^h  = u^h$ and \eqref{fem:tempeq} with $\psi^h = 1$ give
conservation of energy:
\begin{equation}\label{eq:energyconv}
	\mathcal E(u^h,\theta^h):= \int_\Om\left(\textstyle\frac{1}{2}|u^h|^2 + \theta^h\right)(x,t) ~dx
	\equiv  \mathcal E(u^h_0,\theta^h_0)\,.
\end{equation}
Finally, \eqref{fem:tempeq} with $\psi^h = \frac{1}{\theta^h}$ and \eqref{fem:contequation} 
with $\phi^h := \frac{K}{\tau^h}$ give
\begin{eqnarray}
	\frac{d}{dt}\int_\Om \log \theta^h + K \log \tau^h~dx \!\!&=&\!\! 
	\int_\Om \frac{\mu|u^h_x|^2}{\tau^h \theta^h}~dx 
	- \frac{1}{h}\sum_{i= 1}^{N-1} G_i(\tau^h)\jump{L'(\theta^h)}_i
	\left\llbracket\frac{1}{\theta^h}\right\rrbracket_i\nonumber \\
	\!\!&=&\!\!\int_\Om \frac{\mu|u^h_x|^2}{\tau^h \theta^h}~dx 
	+ \frac{1}{h}\sum_{i= 1}^{N-1} \frac{G_i(\tau^h)
	\kappa(\theta^h_{i\dagger})}{(\theta^h_{i*})^2}\jump{\theta^h}_i^2\label{entrp}
\end{eqnarray}
where, for each $i= 1, \ldots, N-1$, $\theta^h_{i\dagger},\, \theta^h_{i*}$ lie between 
$\theta^h|_{E_i}$ and $\theta^h|_{E_{i+1}}$ (obtained from the mean value theorem 
applied to $L'(\theta)$ and $\textstyle\frac{1}{\theta}$, respectively).
Thus the discrete entropy 
\begin{equation*}
	\mathcal {S}(\tau^h,u^h,\theta^h) := \textstyle\frac{1}{2}|u^h|^2 + \theta^h 
	+ K\tau^h  - \log(\theta^h) - K\log(\tau^h)
\end{equation*}
satisfies	
\begin{equation}\label{eq:entropy}
	\frac{d}{dt}\int_\Om \mathcal {S}(\tau^h,u^h,\theta^h)
	~dx +\int_\Om \frac{\mu|u^h_x|^2}{\tau^h \theta^h}~dx 
	+\frac{1}{h}\sum_{i= 1}^{N-1} \frac{G_i(\tau^h)
	\kappa(\theta^h_{i\dagger})}{(\theta^h_{i*})^2}\jump{\theta^h}_i^2= 0\,.
\end{equation}

\begin{Lemma}[Well-posedness of approximation scheme]
Under the same assumptions as in Definition \ref{def:fem} there exists a unique solution $(\tau^h, u^h, \theta^h)$
of the approximation scheme \eqref{fem:contequation}-\eqref{fem:tempeq}. The solution is
defined for all times $t>0$ and there is a constant $C_0$, depending on $h$ and 
$(\tau_0, u_0, \theta_0)$, but independent of time, such that
$$
	C_0^{-1} \leq \tau^h(x,t), \, \theta^h(x,t) \leq C_0 \qquad \forall x \in \Om\,,~ \forall t\geq 0\,.
$$
\end{Lemma}
\bProof
We define the finite element space 
$\vc{S}_h(\Om) := Q_h(\Om)\times V^0_h(\Om) \times Q_h(\Om)$ of $(\tau^h, u^h, \theta^h)$
triples and equip it with the standard $\vc{L}^2(\Om)  = L^2(\Om)\times L^2(\Om) \times L^2(\Om)$ norm. 
Let $\vc{S}_h^+(\Om)$ denote the subset $\vc{S}_h(\Om)\cap \{\,\tau^h, \, \theta^h\, >0\, \}$.
For each $\vc{z}=(z_1,z_2,z_3)\in \vc{S}_h^+(\Om)$ we define $\vc{F}(\vc{z})$ to be the unique element in 
$\vc{S}_h(\Om)$ satisfying
\begin{eqnarray}
	\langle \vc{F}(\vc{z}), \vc{q}\rangle \!\!\!&=&\!\!\! \int_\Om (z_2)_x q_1~dx + 
	\int_\Om \left[\frac{\mu(z_3) (z_2)_x}{z_1} - p(z_3, z_1)\right]\big[(z_2)_x q_3 - (q_2)_x \big]~dx \nonumber\\
	&& \qquad - \frac{1}{h}\sum_{i =1}^{N-1}G_i(z_1)\jump{L'(z_3)}_i\jump{q_3}_i, \label{F} 
\end{eqnarray}
for all $\vc{q} = (q_1,q_2,q_3) \in \vc{S}_h(\Om)$, where $\langle \cdot, \cdot \rangle$ denotes 
the standard inner product on $\vc{L}^2(\Om)$.
Linearity with respect to $\vc{q}$ in \eqref{F} implies that the approximation scheme 
\eqref{fem:contequation}-\eqref{fem:tempeq} corresponds to the ODE
\begin{equation}\label{eq:wellposed-1}
	\vc{z}_t = \vc{F}(\vc{z}), \quad \vc{z}(0) =\vc{z}_0^h = (\tau^h_0, u^h_0, \theta^h_0).
\end{equation}
For the given initial data $\vc{z}_0^h$ we define the compact subset 
\begin{eqnarray*}
	{\vc{K}}(\vc{z}_0^h)\nquad&:=&\nquad \Big\{\vc{z} \in \vc{S}_h^+(\Om);
	\int_\Om z_1~dx \leq \int_\Om \tau^h_0~dx,\, 
	\int_\Om \mathcal{E}(z_2,z_3)~dx \leq \int_\Om \mathcal{E}(u^h_0,\theta^h_0)~dx,\nonumber \\
	&&\quad\qquad\qquad \int_\Om \log \big(z_3 z_1^K \big)~dx 
	\geq \int_\Om \log\big(\theta^h_0 (\tau^h_0)^K\big)~dx \Big\}.
\end{eqnarray*}
It follows from the estimates \eqref{eq:massconv}-\eqref{eq:entropy} that 
any solution $\vc{z}(t)$ of the initial value problem \eqref{eq:wellposed-1} 
belongs to ${\vc{K}}(\vc{z}_0^h)$ for all times it is defined. In turn this gives a time independent 
(but $h$-dependent) upper bound for $\|\vc{z}(t)\|_{\vc L^\infty(\Om)}$, as well as for 
$\|z_1^{-1}(t)\|_{L^\infty(\Om)}$ and $\|z_3^{-1}(t)\|_{L^\infty(\Om)}$.
A direct calculation shows that, thanks to these bounds, $\vc{F}$ is uniformly
Lipschitz continuous on ${\vc{K}}(\vc{z}_0^h)$: there is a constant $L$ (depending on $h$ and
$\vc{z}_0^h$) such that
\[\|\vc{F}(\vc{z}) - \vc{F}(\vc{y})\|_{\vc{L}^2(\Om)} 
\leq L\|\vc{z} - \vc{y}\|_{\vc{L}^2(\Om)}, \quad \forall \vc{z}, \vc{y} \in {\vc{K}}(\vc{z}_0^h)\,.\]
Standard ODE theory now gives existence and uniqueness of the solution $\vc{z}(t)$ to 
\eqref{eq:wellposed-1} for all times $t\geq 0$. I.e., $(\tau^h, u^h, \theta^h)(t)= (z_1(t),z_2(t),z_3(t))$ is the
unique solution to the approximation scheme \eqref{fem:contequation}-\eqref{fem:tempeq} with initial data
$\vc{z}_0^h$.
\qed

\section{Convergence of the approximation scheme}\label{convergence}
In this section we formulate certain (strong) apriori estimates on solutions of the scheme 
\eqref{fem:contequation}-\eqref{fem:tempeq}, and we show that these are sufficient to 
conclude convergence to a weak solution of the system \eqref{eq:contequation}, 
\eqref{eq:momentumeq},  \eqref{eq:tempeq}. As noted in the introduction, this part of the 
analysis works for more general powers $\alpha$ and $\beta$ (see \eqref{mu_kappa}) 
than those assumed in the statement of Theorem \ref{theorem:main}. 
On the other hand we recall that we restrict ourselves to ideal polytropic gases.

\begin{Notation}
Given a set of elements $\{v^h\}$ in a normed space $(X, \|\cdot\|)$, we write ``$v^h\inb X$'' 
to mean that $\|v^h\|$ is bounded independently of $h$. Weak convergence is denoted by $\weak$
and weak-* convergence is denoted $\weakstar$. An over bar denotes weak $L^1$-limit, see e.g.\
\eqref{eq:eff-claim}. Subsequences are not relabeled. 
A zero subscript denotes evaluation at time $t=0$.
\end{Notation}
For later reference we recall the following general results:

\begin{Lemma}\label{lemma:dunford}
Let $\{f^h\}$, $\{g^h\}$ be sequences of functions on a measurable subset $\mathcal{O}$ of  
$\Om \times (0,T)$. Assume that $f^h \inb L^\infty(\mathcal{O})$ with $f^h \rightarrow f$ a.e.\ in 
$\mathcal{O}$, and  $g^h \inb L^1(\mathcal{O})$ with $g^h \weak g$ in $L^1(\mathcal{O})$.
Then $f^hg^h \weak fg$ in $L^1(\mathcal{O})$.
\end{Lemma}
%

\begin{Lemma}\label{lemma:convexity} (\cite[Theorem 2.11]{feireislbook})
Let $O\subset \mathbb{R}^M$, $M\ge 1$, be bounded and open and
let $g\colon \mathbb{R}\to \mathbb{R}$ be continuous and 
convex. Let $\Set{v_n}_{n\ge 1}$ be a sequence in $L^1(O)$ such that $v_n\weak v$ in $L^1(O)$, $g(v_n)\in L^1(O)$ for each 
$n$, and $g(v_n)\weak \overline{g(v)}$ in $L^1(O)$. Then 
$g(v)\le \overline{g(v)}$ a.e.~on $O$ and $g(v)\in L^1(O)$. 
If, in addition, $g$ is strictly convex on an open interval
$(a,b)\subset \mathbb{R}$ and $g(v)=\overline{g(v)}$ a.e.~on $O$, 
then, passing to a subsequence if necessary, 
$v_n(x)\to v(x)$ for a.e.~$x\in \Set{y\in O\mid v(y)\in (a,b)}$.
\end{Lemma}

\begin{Lemma}\label{lemma:aubin_lions} (\cite[Corollary 4]{simon})
Let $X \subset B \subset Y$ be  Banach spaces with
 $X \subset B$ compactly. Then, for $1 \leq p < \infty$, 
$\Set{v:~v \in L^p(0,T;X), v_t \in L^1(0,T;Y)}$ is compactly embedded in $L^p(0,T;B)$.
In the case $p=\infty$, for any $r > 1$, 
$\Set{v:~v \in L^\infty(0,T;X), v_t \in L^r(0,T;Y)}$ is compactly embedded in $C(0,T;B)$.
\end{Lemma}

The main result in this section is:
\begin{Theorem}\label{theorem:convergence}
Assume the constitutive relations \eqref{ideal_polytr} and \eqref{mu_kappa} with $0\leq \beta<2$.
Let $\{(\tau^h, u^h, \theta^h)\}$ be a sequence of functions constructed according
to Definition \ref{def:fem} for $t\in[0,T]$, and assume there exists $C>0$, 
independent of $h$, such that
\begin{eqnarray*}
&{\bf [C1]}:&\quad \left\{\begin{array}{r c l}
		C^{-1} &\leq& \tau^h(x,t), ~\theta^h(x,t) \\
		\mu(\theta^h(x,t)) &\leq& C
	\end{array}\right\}\quad \forall \, (x,t) \in \Om \times [0,T)\\
&{\bf [C2]}:&\quad   \sup_{t \in (0,T)}\|\theta^h(t)\|_{L^2(\Om)} 
		+ \sum_{i=1}^{N-1}\int_0^T \frac{\jump{\theta^h(t)}^2_i}{h}~dt \leq C\\
&{\bf [C3]}:&\quad  u^h \inb L^\infty(0,T;L^2(\Om))\cap L^2(0,T;W^{1,2}_0(\Om))\,.
\end{eqnarray*}
Then, passing to a subsequence if necessary, $(\tau^h, u^h, \theta^h) \rightarrow (\tau, u, \theta)$ a.e 
in $\Om \times (0,T)$, where $(\tau, u, \theta)$ is a weak solution to the compressible Navier--Stokes
system \eqref{eq:contequation} - \eqref{eq:tempeq} in the sense of Definition \ref{def:weak}.
\end{Theorem}
This theorem will be a consequence of Lemmas \ref{lemma:uconv} - \ref{lemma:uconv2}.

\begin{Remark}
Clearly, \Cn{1} is satisfied for $\alpha=0$ while it is equivalent to a uniform upper bound on $\theta^h$ when $\alpha>0$.
\end{Remark}

\noindent
By $\Cn{1}-\Cn{3}$ there
exist functions $\tau$, $u$, $\theta$ such that 
\begin{align}
		u^h &\weakstar u \qquad \textrm{ in }\quad L^\infty(0,T;L^2(\Om)) \cap L^2(0,T;W^{1,2}_0(\Om)) \label{eq:uconv2}\\
		\tau^h &\weak \tau\qquad \textrm{ in }\quad W^{1,2}(0,T;L^2(\Om)) \label{eq:tauconv}\\
		\theta^h &\weakstar \theta\qquad \textrm{ in }\quad L^\infty(0,T;L^{2}(\Om))\,,\label{eq:thetaconv}
\end{align}
if necessary by passing to a subsequence.
Here $\eqref{eq:uconv2}$ follows from \Cn{3}, \eqref{eq:tauconv}
follows from $\tau^h_t = u_x^h$ and \Cn{3},
and \eqref{eq:thetaconv} follows from \Cn{2}.

We will show that the triple $(\tau, u, \theta)$ is a weak solution of the Navier-Stokes system. 
In Section \ref{sssec:strongtemp} we prove that the velocity and temperature converge strongly
by utilizing Lemma \ref{lemma:aubin_lions}.
Then a convexity argument is used repeatedly in Section \ref{sssec:strongtau} to show
that the specific volume converges strongly. At this point it follows that $(\tau, u, \theta)$ is a weak 
solution of \eqref{eq:contequation} and \eqref{eq:momentumeq}. Exploiting this and the fact that the pressure is an 
$L^2$ function, we find that $\{u^h_x\}$ converges strongly. It then follows that the limits
satisfy the temperature equation as well, and Theorem \ref{theorem:convergence} follows.

\subsection{Strong convergence of $u^h$ and $\theta^h$} \label{sssec:strongtemp}
\begin{Lemma}\label{lemma:uconv}
Let $\Set{\tau^h, u^h, \theta^h}$ be as in Theorem \ref{theorem:convergence}. Then
\begin{equation}\label{eq:uconv-time}
	u^h_t \inb L^2(0,T;W^{-1,2}(\Om))\,,
\end{equation}
and consequently, if necessary passing to a subsequence,
\begin{equation}\label{eq:uconv}
	u^h \rightarrow u\qquad  \textrm{ in }\quad L^2(0,T;L^2(\Om))\,.
\end{equation}
\end{Lemma}
\bProof
Fix $v \in L^2(0,T;W^{1,2}_0(\Om))$ and define functions $v^h\in V^0_h(\Om)$ by requiring 
$$
	\int_0^T \int_\Om v^h w^h~dxdt = \int_0^T \int_\Om v w^h~dxdt\qquad \forall w^h \in V^0_h(\Om)\,,
$$
i.e., $v^h$ is the $L^2$-projection of $v$ onto $V^0_h(\Om)$. Then, by \cite{bramble}, 
there exists a constant $C>0$, independent of $h$, such that
\begin{equation*}
	\|v^h\|_{L^2(0,T;W^{1,2}(\Om))} \leq C\|v\|_{L^2(0,T;W^{1,2}(\Om))}.
\end{equation*}
Using $v^h$ as test function in the momentum scheme \eqref{fem:momentumeq} gives
\begin{equation*}
	\begin{split}
		&\left|\int_0^T\int_\Om u^h_t v~dxdt\right|
		 = \left|\int_0^T\int_\Om u^h_t v^h~dxdt \right| \\
		&\quad \leq \|v^h\|_{L^2(0,T;W^{1,2}(\Om))}
			\left[\left\|\frac{\mu(\theta^h)}{\tau^h}u_x^h\right\|_{L^2(0,T;L^2(\Om))} 
			+ \|p(\theta^h, \tau^h)\|_{L^2(0,T;L^2(\Om))}\right] \\
		&\quad\leq C\|v\|_{L^2(0,T;W^{1,2}(\Om))},
	\end{split}
\end{equation*}
where we have used \Cn{1} - \Cn{3}. As $v$ is arbitrary, we conclude  \eqref{eq:uconv-time}.

In the current situation where both \Cn{3} and \eqref{eq:uconv-time} holds, 
Lemma \ref{lemma:aubin_lions}, with $X = W_0^{1,2}(\Om)$, $B=  L^2(\Om)$, and $Y=W^{-1,2}(\Om)$,
yields \eqref{eq:uconv}.
\qed

To analyze the convergence of $\theta^h$ we define a dual mesh $\mathbb{E}_h^{\perp}$ 
with vertices
$$
x_i^\perp = h\big(i - \textstyle\frac{1}{2}\big),~ i= 1, \ldots, N-1, \quad x^\perp_0 = 0, \quad x^\perp_N = L.
$$
We define finite element spaces
$V_h^\perp(\Om)$, $Q_h^\perp(\Om)$ and corresponding projections 
$\Pi_h^{V^\perp}$ and $\Pi_h^{Q^\perp}$ exactly as we did for
$\mathbb{E}_h$. To shorten notation  we write
\begin{equation}\label{eq:dualprojection}
 \widetilde{\phi^h} = \Pi_h^{V^\perp}[\phi]\,.
\end{equation}
\begin{Lemma}\label{lemma:thetaconv}
	Let $\Set{(\tau^h, u^h, \theta^h)}$ be as in Theorem \ref{theorem:convergence}.
	Then, 
	\begin{equation}\label{eq:conv-thetah1}
		\widetilde{\theta^h} \inb L^2(0,T;W^{1,2}(\Om))\,,\qquad \theta^h_t \inb L^1(0,T;W^{-1,1}(\Om))\,,
	\end{equation}
	and $\theta^h \rightarrow \theta$ in $L^p(0,T;L^p(\Om))$, $\forall p\in[1,4)$, if necessary passing to a subsequence.
\end{Lemma}
\bProof
The dual mesh is defined such that:
\begin{equation*}
	h^{-1}\jump{\theta^h}_i = \widetilde{\theta^h}_x, \qquad x \in (x_i^\perp, x_{i+1}^\perp), \quad i=1, \ldots, N-1.
\end{equation*}
and by \Cn{2} there is a constant $C>0$, independent of $h$, such that
\begin{equation}\label{C2}
	\|\widetilde{\theta^h}_x\|_{L^2(0,T;L^2(\Om))}^2 = 
	\int_0^T\sum_{i=1}^{N-1}\frac{\jump{\theta^h}_i^2}{h}~dt \leq C\,.
\end{equation}
Since $\|\widetilde{\theta^h}(t) - \theta^h(t)\|_{L^2(\Om)} \leq Ch\|\widetilde{\theta_x^h}(t)\|_{L^2(\Om)}$, we thus have
\begin{equation*}
	\|\widetilde{\theta^h} - \theta^h\|_{L^2(0,T;L^2(\Om))} \leq Ch\|\widetilde{\theta^h}_x\|_{L^2(0,T;L^2(\Om))} \leq Ch\,,
\end{equation*}
and \Cn{2} gives
\begin{equation*}
	\|\widetilde{\theta^h}\|_{L^2(0,T;L^2(\Om))} 
	\leq \|\widetilde{\theta^h} - \theta^h\|_{L^2(0,T;L^2(\Om))}
	+ \|\theta^h\|_{L^2(0,T;L^2(\Om))}  \leq C(1+h)\,,
\end{equation*}
which proves $\eqref{eq:conv-thetah1}_1$. For later reference we note that,
again by \eqref{C2}, $\widetilde{\theta^h} \inb L^2(0,T;L^\infty(\Om))$, whence the same 
holds for $\theta^h$. Together with \Cn{2} and the Cauchy-Schwarz inequality this gives
\begin{equation}\label{eq:interpolation}
	\theta^h \inb L^4(0,T;L^4(\Om))\,,
\end{equation}
which in turn shows that $\kappa(\theta^h) \inb L^2(0,T;L^2(\Om))$, since $\beta < 2$.

For $\eqref{eq:conv-thetah1}_2$ we fix $\psi \in W^{1,\infty}(\Om)$ and define the discrete effective
viscous flux by $\Fh := \frac{\mu(\theta^h)}{\tau^h}u_x^h - p(\theta^h, \tau^h)$.
By \Cn{1} - \Cn{3}, and the Cauchy-Schwarz inequality we have
$\Fh  \inb L^2(0,T;L^2(\Om))$.
Using $\Pi_h^Q\psi$ as test function
in the temperature approximation \eqref{fem:tempeq}, we deduce
\begin{equation}\label{eq:thetatime}
	\begin{split}
		&\Big|\int_\Om \theta^h_t \psi~dx\Big| = \Big|\int_\Om \Fh u^h_x \psi~dx
		 - \frac{1}{h}\sum_{i =1}^{N-1}G_i(\tau^h)\jump{L'(\theta^h)}_i\jump{\Pi_h^Q \psi}_i\Big| \\
		& \leq \|\psi\|_{L^\infty(\Om)}\|\Fh u^h_x\|_{L^1(\Om)} + \|\psi_x\|_{L^\infty(\Om)}\sum_{i=1}^{N-1}
		G_i(\tau^h)\big|\jump{L'(\theta^h)}_i\big| \\
		&\leq \|\psi\|_{W^{1,\infty}(\Om)}\Big[\|\Fh\|_{L^2(\Om)}\|u_x^h\|_{L^2(\Om)}
		 + C\|\kappa(\theta^h)\|_{L^2(\Om)}\Big[\sum_{i=1}^{N-1}\frac{1}{h}\jump{\theta^h}_i^2\Big]^\frac{1}{2}\Big],
	\end{split}
\end{equation}
where we have used \Cn{1} and that $L''(z) = \kappa(z)$. 
Integrating \eqref{eq:thetatime} and applying \Cn{2}-\Cn{3}, together with the bounds on $\Fh$ and 
$\kappa(\theta^h)$, then gives $\eqref{eq:conv-thetah1}_2$.

Next we establish compactness of $\{\theta^h\}$ with respect to the spatial variable.
Let $\xi \in (-L,L)$ be arbitrary and define the domain 
$\Om_\xi = (|\xi|,L-|\xi|)$.
If $h \leq |\xi|$, fix any $x \in \Om_\xi$, let $E_i$ be the element containing $x$, and
let $E_j$ be the element containing $x - \xi$. Then
\begin{equation*}
	\begin{split}
	&|\theta^h(x,t) - \theta^h(x-\xi,t)|^2 
	= |\theta^h(t)|_{E_i} - \theta^h(t)|_{E_j}|^2 \\
	&= \left|\pm \sum_{k=i}^{j-1}h^\frac{1}{2}h^{-\frac{1}{2}}\jump{\theta^h}_k\right|^2 
	\leq \left[\sum_{k=i}^{j-1}h \right]\left[\sum_{k= 1}^{N-1}\frac{\jump{\theta^h}_k^2}{h} \right] 
	\leq C(|\xi|+ h)\left[\sum_{k= 1}^{N-1}\frac{\jump{\theta^h}_k^2}{h} \right]\,.
\end{split}
\end{equation*}
Integrating over $\Om_\xi\times (0,T)$, and using $h \leq |\xi|$ together with \Cn{2}, we obtain
\begin{equation}\label{eq:translationwithh}
	\|\theta^h(x,t) - \theta^h(x-\xi, t)\|_{L^2(0,T;L^2(\Om_\xi))} 
	\leq C|\xi|^\frac{1}{2}\,.
\end{equation}
On the other hand, if $|\xi| \leq h$, then for each $0\leq t\leq T$ we have (for $h<1$)
$$
\int_0^T\!\!\int_{\Om_\xi}\!\! |\theta^h(x,t) - \theta^h(x-\xi,t)|^2dxdt 
= \int_0^T\!\! |\xi| \sum_{i=1}^{N-1}\jump{\theta^h}_i^2\,dt
\leq |\xi| \int_0^T\!\sum_{i=1}^{N-1}\frac{\jump{\theta^h}_i^2}{h}\, dt\,,
$$
whence \Cn{2} shows that \eqref{eq:translationwithh} holds also for $|\xi| \leq h$. 
From the translation estimate \eqref{eq:translationwithh}, it follows that 
 $\theta^h \inb L^2(0,T;W^{\tfrac{1}{2} -\epsilon,2}(\Om))$, see \cite{Adams}.

Consequently, Lemma \ref{lemma:aubin_lions}, with $X= W^{\tfrac{1}{2}- \epsilon, 2}(\Om)$, $\epsilon >0$, 
$B= L^2(\Om)$, and $Y= W^{-1,1}(\Om)$, yields
\begin{equation}\label{theta-L2}
	\theta^h \rightarrow \theta \qquad \textrm{in $L^2(0,T;L^2(\Om))$ }\,.
\end{equation}
Finally,  \eqref{eq:interpolation} and \eqref{theta-L2} 
concludes the proof.\qed

\noindent
We observe that \eqref{theta-L2} and \Cn{1} yields:
\begin{equation}\label{eq:conv1}
	\mu(\theta^h) \rightarrow \mu(\theta) \qquad \text{in}\quad L^p(\Om\times (0,T)) \qquad \forall\, 1\leq p< \infty.
\end{equation}

\subsection{Strong convergence of $\tau^h$}\label{sssec:strongtau}
We first note that $\tau^h_t = u^h_x$ a.e.\ by \eqref{fem:contequation}, whence
\eqref{eq:uconv2} and \eqref{eq:tauconv} show that  $\tau,\, u$ satisfy \eqref{weak:masseq}. 
Passing to the limit a.e. in the other equations requires strong convergence of $\tau^h$. We
begin by observing that, due to \Cn{1}, we may divide by $\tau$ in 
\eqref{weak:masseq} to obtain the ''renormalized'' equation (used below in \eqref{eq:thepoint})
\begin{equation}\label{eq:renorm}
	(\log \tau)_t = \frac{u_x}{\tau} \qquad \textrm{a.e.\ in }\Om \times (0,T)\, .
\end{equation}
We proceed to adapt a simplified version of the convexity arguments in \cite{lions}, \cite{f} 
to the Lagrangian setting and obtain the required strong convergence.

\begin{Lemma}\label{lemma:tauconv}
Let $\Set{\tau^h, u^h, \theta^h}$ be as in Theorem \ref{theorem:convergence}.
Then,
\begin{equation*}
	\tau^h \rightarrow \tau \qquad \textrm{a.e.\ in } \Om \times (0,T)\,.
\end{equation*}
\end{Lemma}
\bProof
Step 1: Recalling the notation for weak $L^1$-limits we claim that 
\begin{equation}\label{eq:eff-claim}
	\overline{\left(\frac{u_x}{\tau}\right)} - \frac{u_x}{\tau}\geq 0 \qquad\textrm{a.e.\ in }\Om \times (0,T)\,.
\end{equation}
For the proof we fix $\psi \in C_c^\infty(\Om\times (0,T))\cap \Set{\psi \geq 0 }$ and set
$$
	v^h(x,t):= \int_0^x \Big[\psi \frac{\tau - \tau^h}{\tau}\Big](y,t)\, dy
	- \frac{x}{|\Omega|}\int_\Omega \Big[\psi \frac{\tau - \tau^h}{\tau}\Big](y,t)\, dy\,,
$$
where the integrand is differentiable in time since
\begin{equation}\label{eq:intdiff}
	\Big[\frac{\tau - \tau^h}{\tau}\Big]_t = \frac{\tau^h u_x - \tau u_x^h}{\tau^2} \inb L^2(0,T;L^1(\Om))\,.
\end{equation}
Hence, $v^h_t \inb L^2(0,T;L^\infty(\Om))$, and by \eqref{eq:uconv2} and \eqref{eq:tauconv} we have
\begin{equation}\label{eq:testfunc-conv2}
	v^h_t \weakstar 0 \qquad \text{in }L^2(0,T;L^\infty(\Om))\,.
\end{equation}
We then use $\Pi^V_h[v^h]$ as test function \eqref{fem:momentumeq}
and integrate in time to obtain 
\begin{eqnarray}\label{eq:eff-1}
	&&\nquad\int_0^T\int_\Om \left[\frac{\mu(\theta^h)u_x^h}{\tau^h} - \frac{\mu(\theta^h)u_x^h}{\tau}\right]\!\psi\, dxdt 
	- \int_0^T \int_\Om \psi \left[p(\theta^h, \tau^h) - p(\theta^h, \tau)\right]dxdt \nonumber \\
	&&\quad = \int_0^T\left[\int_\Om \frac{\mu(\theta^h)u^h_x}{\tau^h} - p(\theta^h, \tau^h)~dx\right]\!\!
				\left[\frac{1}{|\Om|}\int_\Om \psi \frac{\tau - \tau^h}{\tau}~dy\right]dt \nonumber\\
	&&\quad + \int_0^T \int_\Om u^h v^h_t ~dxdt 
	+ \int_0^T \int_\Om u^h \left[\Pi^V_h[v^h_t] - v^h_t \right]dxdt 
=: \sum_{i=1}^3 I_i\,, 
\end{eqnarray}
where we have used \eqref{proj1} and \eqref{proj2}.
We observe that $I_1, I_2 \rightarrow 0$ as $h \rightarrow 0$, due to \eqref{eq:testfunc-conv2} and 
Lemma \ref{lemma:uconv}. The $I_3$ term satisfies
\begin{equation*}
	|I_3| \leq \|u^h\|_{L^2(0,T;L^\infty(\Om))}\big\|\Pi^V_h[v^h_t] - v^h_t\big\|_{L^2(0,T;L^1(\Om))},
\end{equation*}
where the first factor is bounded by \Cn{3}. A standard interpolation estimate \cite{bs}, together with \eqref{eq:intdiff}, gives
\begin{equation*}
	\big\|\Pi^V_h[v^h_t] - v^h_t\big\|_{L^2(0,T;L^1(\Om))} \leq Ch\|v_{tx}^h\|_{L^2(0,T;L^1(\Om))} \leq Ch\,.
\end{equation*}
Recalling \eqref{eq:conv1} and Lemma \ref{lemma:thetaconv}, and sending $h \rightarrow 0$ in \eqref{eq:eff-1} we get
\begin{equation*}
	\lim_{h \rightarrow 0}  \int_0^T \int_\Om \mu(\theta)\Big[\frac{u^h_x}{\tau^h} - \frac{u_x}{\tau}\Big]\psi\, dxdt \\
	= \lim_{h \rightarrow 0} \int_0^T\int_\Om K\theta \Big[\frac{1}{\tau^h} - \frac{1}{\tau}\Big]\psi\, dxdt
	\geq 0,
\end{equation*}
where the last inequality follows from Lemma \ref{lemma:convexity} and the convexity of $z \mapsto \frac{1}{z}$.
We thus have
\begin{equation*}
	\mu(\theta)\left[\overline{\left[\frac{u_x}{\tau}\right]} - \frac{u_x}{\tau}\right] \geq 0\quad \textrm{ a.e in }\Om\times (0,T)\,,
\end{equation*}
and the claim \eqref{eq:eff-claim} follows from the lower bound \Cn{1} on $\mu(\theta)$.

\noindent
Step 2: By \eqref{eq:renorm} and \eqref{eq:eff-claim} it follows that
\begin{equation}\label{eq:thepoint}
	\begin{split}	
		\big[\overline{\log \tau} - \log \tau\big]_t = 
		\overline{\left[\frac{u_x}{\tau}\right]} - \frac{u_x}{\tau} \geq 0 \qquad \textrm{a.e.\ in }(0,T) \times \Om\, .
	\end{split}
\end{equation}
As $\log \tau^h_0 \rightarrow \log \tau_0$ a.e.\ in $\Om$ we conclude that
\begin{equation*}
	\overline{\log \tau}(x,t) - \log \tau(x,t) \geq 0\qquad \textrm{ for a.e } (x,t) \in \Om \times (0,T)\,.
\end{equation*}
On the other hand, Lemma \ref{lemma:convexity} shows that $\overline{\log \tau}(x,t) - \log \tau(x,t) \leq 0$ a.e.\ 
Thus,  $\overline{\log \tau} = \log \tau$ a.e.\ in $ \Om\times(0,T)$, and the conclusion follows by an application 
of the last part of Lemma \ref{lemma:convexity}.
\qed

\subsection{Concluding the proof of Theorem \ref{theorem:convergence}}
In view of the two previous lemmas and \Cn{1},
\begin{equation}\label{eq:conv2}
	\begin{split}
		p(\theta^h, \tau^h) \rightarrow p(\theta, \tau) \qquad &\text{in } L^2(0,T;L^2(\Om)), \\
		\frac{\mu(\theta^h)}{\tau^h} \rightarrow \frac{\mu(\theta)}{\tau} \qquad &\text{in } L^p(\Om\times(0,T)),\, \forall 1\leq p< \infty\, .
	\end{split}
\end{equation}
Given $\phi \in C^\infty_c(\Om\times[0,T))$ we use \eqref{fem:momentumeq} with 
$v^h=\pivh{\phi}$ and integrate in  time:
\begin{equation}\label{eq:mom-conv1}
	\int_0^T\int_\Om u^h\pivh{\phi_t}~dxdt + \int_\Om u_0^h \pivh{\phi_0}~dx
	= \int_0^T\int_\Om \mathcal F^h\piqh{\phi_x}~dxdt\,.
\end{equation}
Due to the regularity of $\phi$, standard interpolation estimates show, together with
\Cn{0}, \eqref{eq:uconv2}, and  \eqref{eq:conv2}, that  the limit $h \rightarrow 0$ in \eqref{eq:mom-conv1}
results in \eqref{weak:momentumeq}. By now we have established the a.e.\ convergence $(\tau^h, u^h, \theta^h)\to 
(\tau, u, \theta)$, and that the latter triple satisfies \eqref{weak:masseq} and \eqref{weak:momentumeq}.
It remains to prove that $(\tau, u, \theta)$ also satisfies \eqref{weak:tempeq}.
In order to pass to the limit in the temperature scheme, the following lemma is essential;
its proof exploits the fact that $(\tau, u, \theta)$ satisfies \eqref{weak:momentumeq}.
\begin{Lemma}\label{lemma:uconv2}
Let $\Set{\tau^h, u^h, \theta^h}$ be as in Theorem \ref{theorem:convergence}.
Then,
\begin{equation*}
	\lim_{h \rightarrow 0}\left[\sup_{t \in (0,T)}\int_\Om \frac{|u^h - u|^2}{2}~dx 
	+ \int_0^T\int_\Om \frac{\mu(\theta^h)}{\tau^h}|u_x^h - u_x|^2~dxdt\right] = 0.
\end{equation*}
\end{Lemma}
\bProof
Using $u^h$ as test function in \eqref{fem:momentumeq} and integrating in time we get
\begin{equation}\label{eq:uxconv1}
\int_\Om \frac{|u^h(t)|^2}{2}~dx = \int_\Om \frac{|u_0^h|^2}{2}~dx 
- \int_0^t\int_\Om \mu(\theta^h) \frac{|u_x^h|^2}{\tau^h} - p(\theta^h, \tau^h)u_x^h\,dxdt\,.
\end{equation}
Sending $h \rightarrow 0$ in \eqref{eq:uxconv1}, and using \eqref{eq:conv2}, gives
\begin{equation*}
	 \lim_{h \rightarrow 0}\left[\frac{\|u^h(t)\|_2^2}{2}
	 +\int_0^t \int_\Om \mu(\theta^h) \frac{|u_x^h|^2}{\tau^h}\, dxdt\right] 
	 =  \frac{\|u_0^h\|_2^2}{2} + \int_0^t \int_\Om p(\theta, \tau)u_x\,dxdt.
\end{equation*}
On the other hand, by Lemma \ref{lemma:uconv}, $u_t \inb L^2(0,T;W^{-1,2}(\Om))$. This and
the weak form of the momentum equation \eqref{weak:momentumeq} yield the weak form
\begin{equation*}
	\int_0^t \left\langle u_t, v \right\rangle_{\langle W^{-1,2},W^{1,2}_0\rangle} ~dt 
	=- \int_0^t\int_\Om \Big[\frac{\mu(\theta)u_x}{\tau} - p(\theta, \tau)\Big]v_x~dxdt,
\end{equation*}
valid for all $v \in L^2(0,T;W^{1,2}_0(\Om))$.
Applying this with $v = u$ gives
\begin{equation*}
	\frac{\|u(t)\|_2^2}{2}+ 	 \int_0^t\int_\Om \mu(\theta) \frac{|u_x|^2}{\tau} ~dxdt
	=  \frac{\|u_0\|_2^2}{2} + \int_0^t\int_\Om p(\theta, \tau)u_x~dxdt\,,
\end{equation*}
whence 
\begin{equation*}\label{eq:turisttaske}
	\begin{split}
	& \lim_{h \rightarrow 0}\left[\frac{\|u(t)\|_2^2}{2} - \frac{\|u_0\|_2^2}{2} 
	+\int_0^t\int_\Om \mu(\theta^h) \frac{|u_x^h|^2}{\tau^h}-\mu(\theta) \frac{|u_x|^2}{\tau} ~dxdt\right] 
	=0\quad\forall t \in (0,T)\,.
	\end{split}
\end{equation*}
Next, $u_x^h \weak u_x$ in $L^2(0,T;L^2(\Om))$ and $\frac{\mu(\theta^h)}{\tau^h} \rightarrow \frac{\mu(\theta)}{\tau}$
in $L^p(0,T;L^p(\Om))$ for any $p < \infty$. 
Consequently, in view of Lemma \ref{lemma:dunford},
$$
\lim_{h\rightarrow 0}\int_0^T\int_\Om \frac{\mu(\theta^h)}{\tau^h}u_x^hu_x~dxdt 
= \int_0^T\int_\Om\frac{\mu(\theta)}{\tau}u_x^2~dxdt\, ,
$$
such that 
\begin{equation*}
	\begin{split}
		 &\lim_{h \rightarrow 0}\left[\sup_{t \in (0,T)}\int_\Om \frac{|u^h - u|^2}{2}\,dx 
		 + \int_0^T\int_\Om \frac{\mu(\theta^h)}{\tau^h}|u_x^h - u_x|^2\, dxdt\right]  \\
		&=  \lim_{h \rightarrow 0}\left[\sup_{t \in (0,T)}\int_\Om \frac{|u^h|^2 - |u|^2}{2}\, dx 
		+\int_0^T\int_\Om \mu(\theta^h) \frac{|u_x^h|^2}{\tau^h}-\mu(\theta) \frac{|u_x|^2}{\tau} \,dxdt\right]
		=0\,.
	\end{split}
\end{equation*}
As the left hand side is non-negative this concludes the proof. 
\qed

\noindent
Finally, fix any $\phi \in C_c^\infty(\Om \times [0,T))$.
A calculation using \eqref{eq:dualprojection} shows that
\begin{align}\label{eq:purevirtue}
	\sum_{i=1}^{N-1}\frac{1}{h}G_i(\tau^h)\jump{L'(\theta^h)}_i\jump{\piqh {\phi}}_i 
	&= \int_\Om \frac{\kappa(\theta_{i\dagger}^h)}{\Pi_h^{Q^\perp}\left[\tau^h\right]} 
	\big(\widetilde{\theta^h}\big)_x \widetilde{\big(\piqh{\phi}\big)}_x ~dx,
\end{align}
where $\theta^h_{i\dagger}$ is as in \eqref{entrp}. Now using $\piqh{ \phi}$ as test function in \eqref{fem:tempeq} we obtain 
\begin{align}\label{eq:finally}
		&\nquad\int_0^T\int_\Om \theta^h\phi_t~dxdt - \int_\Om \theta^h_0 \phi_0\,dx 
		- \int_0^T\int_\Om \frac{\kappa(\theta_{i\dagger}^h(x))}{\tau^h} (\widetilde{\theta^h})_x \widetilde{(\piqh{\phi})}_x\,dxdt \nonumber\\
		&\nquad= \int_0^T\int_\Om \Big[p(\theta^h, \tau^h)u^h_x-\frac{\mu(\theta^h)|u_x|^2}{\tau^h} 
		-\frac{\mu(\theta^h)\big[|u^h_x|^2 - |u_x|^2\big]}{\tau^h}\Big]\phi\,dxdt.
\end{align}
On each interval $(x_i^\perp, x_{i+1}^\perp)$, $\kappa(\theta^h_{i\dagger})$ is a linear combination 
of $\kappa(\theta^h_{i})$ and $\kappa(\theta^h_{i+1})$. Hence, by Lemma \ref{lemma:thetaconv}
and since $\beta<2$, $\kappa(\theta^h_{i\dagger}) \rightarrow \kappa(\theta)$ 
in $L^2(0,T;L^2(\Om))$. Thus $\kappa(\theta_{i\dagger}^h(x)) (\widetilde{\theta^h})_x \weak \kappa(\theta)\theta_x$
in $L^1(0,T;L^1(\Om))$. Also, as
$\frac{1}{\tau^h} \rightarrow \frac{1}{\tau}$ in $L^p(0,T;L^p(\Om))$ $\forall p\in[1, \infty)$,
Lemma \ref{lemma:dunford} gives
$$
\frac{\kappa(\theta_{i\dagger}^h)}{\tau^h} (\widetilde{\theta^h})_x \weak \frac{\kappa(\theta)}{\tau}\theta_x 
\qquad \text{ in $L^1(0,T;L^1(\Om))$}\,.
$$
Letting $h \rightarrow 0$ in \eqref{eq:finally} and using Lemma \ref{lemma:uconv2} and \eqref{eq:conv2}, 
results in \eqref{weak:tempeq}. This concludes the proof of Theorem \ref{theorem:convergence}.


%
%


\section{Energy bounds and proof of Theorem \ref{theorem:main}}\label{sec:proof}
We now restrict to constitutive relations as described in Theorem \ref{theorem:main}, for which the pointwise 
estimates \eqref{main:pointwise} are established in the next section. Taking these for granted
we now prove the integral bounds \eqref{main:gradient}. As demonstrated below these
suffice, via Theorem \ref{theorem:convergence}, to establish Theorem \ref{theorem:main}.
For reference we note:
\begin{description}
\item [{\bf [A1]}: ] We have $\alpha=0$ and $\beta\in [0,\textstyle\frac{3}{2})$ such that $\Fh$ and $L$ are given by
	\[\Fh = \frac{\bar\mu}{\tau^h}u_x^h - p(\theta^h, \tau^h)\qquad\mbox{and}\qquad 
	L(\theta)=\frac{\bar \kappa \, \theta^{\beta+2}}{(\beta+1)(\beta+2)}\,;\] 
\item[{\bf [A2]}: ] $(\tau^h, u^h, \theta^h)$, $h>0$, are solutions of the scheme \eqref{fem:contequation}-\eqref{fem:tempeq} 
	with initial data $(\tau_0^h, u_0^h, \theta_0^h)$ given by \eqref{data_h} with $(\tau_0, u_0, \theta_0)$ as in
	Theorem \ref{theorem:main}.
\end{description}
For a fixed $T\in(0,\infty)$ we let $C$, $\tilde C$, etc.\ be numbers that depend on $T$, system parameters 
($\Omega$, $\bar \mu$, $K$), and the initial data, but that are independent of $h$.
\begin{Lemma}[Pointwise estimates]\label{lemma:pointwise}
Assume {\bf [A1]} and {\bf [A2]}. Then there exists a number $C>0$ which is independent of $h$ and such that:
\begin{equation}\label{tau_bnd}
	C^{-1}\leq \tau^h(x,t) \leq C, \quad \forall (x,t) \in \Om\times (0,T),
\end{equation}
\begin{equation}\label{theta_bnd}
	C^{-1} \leq \theta^h(x,t), \quad \forall (x,t) \in \Om\times (0,T),
\end{equation}
\begin{equation}\label{theta_int_bnd}
		\int_0^T \|\theta^h(t)\|_{L^\infty(\Om)}~dt \leq C\,.
\end{equation}
\end{Lemma}
The proof of \eqref{tau_bnd}-\eqref{theta_int_bnd} is essentially the same as in \cite{ks}. 
Minor adjustments are required to treat the particular scheme \eqref{fem:contequation}-\eqref{fem:tempeq}
and to incorporate $\theta$-dependence in the heat conductivity $\kappa$. 
For completeness we include the proof in Section \ref{sec:pointwise}. 
We note that these pointwise bounds do not seem to generalize in any simple way 
to the case of $\theta$-dependent viscosities.

We proceed to state the discrete analogue of \eqref{main:gradient} (Lemma \ref{lemma:gradient} below),
and then show how this is used together with Lemma \ref{lemma:pointwise} to prove Theorem \ref{theorem:main}. 
The remaining parts of this section detail the proof of Lemma \ref{lemma:gradient}. A few technical lemmas 
are collected at the end of the section.

\begin{Lemma}[Energy estimates]\label{lemma:gradient}
Define
\begin{equation*}
	\begin{split}
	\Ah(t) &:= \sup_{s \in (0,t)}\int_\Om  |u^h_x|^2~dx + \int_0^t\int_\Om |u^h_t|^2~dxds, \\
	\Bk(t) &:= \sup_{s \in (0,t)}\int_\Om  L(\theta^h) ~dx 
	+ \frac{1}{h}\sum_{i=1}^{N-1} \int_0^t  G_i(\tau^h) \jump{L'(\theta^h)}_i^2 ~ds, \\
	\Dk(t) &:= \sup_{s \in (0,t)}\frac{1}{2h}\sum_{i=1}^{N-1}  G_i(\tau^h)\jump{L'(\theta^h)}_i^2
	+\int_0^t \int_\Om \kappa(\theta^h)|\theta_t^h|^2~dxds  \,.
	\end{split}
\end{equation*}	
Then there is a number $C >0$, independent of $h$,  such that
\begin{equation}\label{main:basic}
	\Ah(t) + \Bk(t) + \Dk(t) \leq C \qquad \forall t \in (0,T)\, . 
\end{equation}
\end{Lemma}

\subsection{Proof of Theorem \ref{theorem:main}}\label{biggie}
We now take Lemmas \ref{lemma:pointwise} and \ref{lemma:gradient} for granted,
and we verify that these are sufficient to verify the conditions in 
Theorem \ref{theorem:convergence} (with $\alpha = 0$). First, our assumptions on 
the initial data in Theorem \ref{theorem:main} are stronger than the corresponding 
conditions in Theorem \ref{theorem:convergence}. Next, $\Cn{1}$ is an immediate 
consequence of Lemma \ref{lemma:pointwise}, \Cn{2} follows from the bound on $\Bh$ 
together with the pointwise estimates of Lemma \ref{lemma:pointwise},
and $\Cn{3}$ follows from the bound on $\Ah$. We can thus apply Theorem \ref{theorem:convergence} 
and conclude the existence and convergence parts of Theorem \ref{theorem:main}. 
At this point, Lemma \ref{lemma:pointwise}  yields \eqref{main:pointwise}.

It remains to verify \eqref{main:gradient}, and we consider the first term in the 
$\mathcal D(t)$-functional in detail. By the $\Dk(t)$-bound in \eqref{main:basic}, we have as in 
\eqref{eq:purevirtue} that
\begin{equation}\label{eq:Kjell-T-Ring}
 	 \int_\Om \frac{|\kappa(\theta_{i\dagger}^h)|^2}{\Pi_h^{Q^\perp}\left[\tau^h\right]} |\widetilde{\theta^h}_x|^2\, dx
 	= \frac{1}{h}\sum_{i=1}^{N-1}  G_i(\tau^h)\jump{L'(\theta^h)}_i^2\leq C\qquad \forall t\in (0,T)\,.
\end{equation}
In view of Lemma \ref{lemma:thetaconv}, this estimate gives
$$
	\widetilde{\theta^h}_x \weakstar \theta_x \qquad \text{in }L^\infty(0,T;L^2(\Om)).
$$ 
Clearly, $\frac{1}{\tau}$, $\theta \in L^\infty(0,T;L^\infty(\Om))$, such that
\begin{equation*}
	\begin{split}
		&\sup_{t \in (0,T)}\int_\Om |L'(\theta)_x|^2~dx = \sup_{t \in (0,T)}\int_\Om \frac{|\kappa(\theta)|^2}{\tau}|\theta_x|^2~dx \\
		&\qquad \leq C\|\theta\|_{L^\infty(0,T;L^\infty)}\Big\|\frac{1}{\tau}\Big\|_{L^\infty(0,T;L^\infty(\Om))}\|\theta_x\|_{L^\infty(0,T;L^2(\Om))} 
		\leq C.
	\end{split}
\end{equation*}
The other terms in \eqref{main:basic} are treated similarly
and hence we conclude \eqref{main:gradient}.\qed

\subsection{Proof of Lemma \ref{lemma:gradient}}\label{sec:gradient}
\begin{Lemma}\label{lemma:basich1}
Assume {\bf [A1]} and {\bf [A2]}. Then there is a number $C>0$, independent of $h$, such that
\begin{equation}\label{energy_0}
	\sup_{t \in (0,T)} \int_\Omega |u^h(x,t)|^2~dx 
	+ \int_0^T \int_\Omega |u^h_x(x,t)|^2~dxdt \,\,\leq \,\,C\,.
\end{equation}
\end{Lemma}
\bProof
Applying \eqref{fem:momentumeq} with $v^{h} = u^{h}$ and $\mu\equiv \bar \mu$, yields
\begin{equation}
\label{kylling}
	\frac{d}{dt}\int_{\Omega}\frac{|u^{h}|^2}{2} ~dx 
	+ \int_{\Omega}\frac{\bar \mu}{\tau^h}|u^h_x|^2 ~dx 
	= \int_{\Omega}p(\tau^{h}, \theta^{h}) u^h_x ~dx\,.
\end{equation}
Applying the Cauchy-Schwarz inequality with a suitable parameter $\epsilon$ we have
\[\int_{\Om} p(\tau^{h}, \theta^{h})u^h_x ~dx  
\leq \epsilon \int_{\Omega}\frac{\bar \mu}{\tau^{h}}|u^h_x|^2 ~dx + 
\frac{\tilde C}{\epsilon}\|\frac{1}{\tau^{h}}\|_{L^\infty(\Omega)}\|\theta^{h}\|_{L^\infty(\Omega)}
\int_{\Omega}\theta^{h} ~dx\,.\]
We choose $\epsilon$ small enough that the first term can be absorbed on the left-hand 
side in \eqref{kylling}. We then apply \eqref{eq:energyconv} together with \eqref{tau_bnd}. 
Integrating in time and using \eqref{theta_int_bnd}, yield \eqref{energy_0} with 
a suitable $C$.
\qed

\subsubsection*{Bound for $\Ah$}
Applying \eqref{fem:momentumeq} with the test function 
$v^h(x,t) := u_t^h(x,t)$, and $\mu\equiv \bar \mu$, gives
\begin{equation*}
	\begin{split}
	\int_\Om |u_t^h|^2~dx 
	&= -\int_\Om \frac{ \bar \mu u_x^h u_{tx}^h}{\tau^h} -  p^h u_{tx}^h~dx \\
	&= - \bar \mu\frac{d}{dt}\left[\int_\Om \frac{  |u^h_x|^2}{2\tau^h}~dx \right]
	+ \int_\Om \frac{\bar \mu  |u_x^h|^2}{2} \left[\frac{1}{\tau^h}\right]_t
	+  p^hu^h_{xt}~dx\, ,
	\end{split}
\end{equation*}
where $p^h:=p(\theta^h,\tau^h)$. Integrating in time and using $\tau^h_t=u^h_x$ we obtain
\begin{eqnarray}
	&&\nquad \nqquad\bar\mu \int_\Om \frac{|u^h_x|^2}{2\tau^h}(t)~dx
	+\int_0^t\int_\Om  |u_t^h|^2~dxds\label{eq:a-basic} \\
	&&\nquad\nqquad= \int_\Om \frac{|u^h_x|^2}{2\tau^h}(0)~dx-
	\int_0^t \frac{\bar \mu  (u_x^h)^3}{2(\tau^h)^2}  + p^h_t u_x^h ~dxds
	+ \int_\Om p^h u^h_x(s) ~dx\Big|_{0}^{t} =: \sum_{i=1}^4 I_i\,,\nonumber
\end{eqnarray}
where we have also applied integration by parts to the pressure term.
From the requirements on the initial data, we have that $I_1 \leq \tilde C$.
Next, by adding and subtracting the positive term $ |u_x^h|^2 p^h/\tau^h$, we have
\[I_2 = - \frac{1}{2}\int_0^t\int_\Om  \frac{\bar \mu  (u_x^h)^3}{(\tau^h)^2}  ~dxds
\leq \left|\int_0^t\int_\Om\frac{ (u_x^h)^2 \Fh}{\tau^h}~dxds\right|\,. \]
Applying Lemma \ref{lemma:fbound}, with $\phi = \frac{|u_x^h|^2}{|\tau^h|}$, gives 
\begin{eqnarray*}
	I_2 &\leq& \tilde C\left\{\epsilon \left[1 +  \Ah(t)\right]
	+ \frac{1}{\epsilon}\int_{0}^t 
	\left[\int_{\Omega}\frac{|u_x^h|^2}{\tau^h}(x,s)~dx\right]^2 ds\right\}
	\nonumber \\
	&\leq& \tilde C\left\{\epsilon\left[1 +  \Ah(t)\right] +
	\frac{1}{\epsilon}\int_0^t \Ah(s)\|u_x^h(s)\|_{L^2(\Om)}^2~ds\right\}\,.\label{eq:a-i2}
\end{eqnarray*}
Next, $p^h=\frac{ K\theta^h}{\tau^h}$, such that
\[I_3 = K\!\!\int_0^t\int_\Om \!\!\frac{ |u^h_x|^2\theta^h}{|\tau^h|^2}
- \frac{  u^h_x\theta^h_t}{\tau^h}~dxds \leq \tilde C\left[\int_0^t \Ah(s)\|\theta^h(s)\|_{L^\infty(\Om)}~ds 
+ \Dk(t)^\frac{1}{2}\right]\,,\]
where we have used \eqref{tau_bnd}, \eqref{theta_bnd}, and Lemma \ref{lemma:basich1}.
To bound $I_4$ we use the Cauchy-Schwarz inequality with a parameter
and \eqref{theta_bnd} together with the requirements on the initial data:
\[I_4 = \left.\int_\Om  p^h u^h_x (s)~dx\right|_{s= 0}^{s=t} 
\leq \tilde C\left[\epsilon \Ah(t) + \frac{1}{\epsilon} \Bk(t)\right] + \|p^h(0)\|_{L^2(\Om)}\|u_x(0)\|_{L^2(\Om)}.\]
Using these bounds in \eqref{eq:a-basic},
taking the supremum over times in $(0,t)$, applying \eqref{tau_bnd}, 
and choosing $\epsilon$ sufficiently small, we obtain
\[\Ah(t) \leq \tilde C\left\{1 + \Bk(t) +\Dk(t)^\frac{1}{2}
+\int_0^t \Ah(t)\left[\|u^h_x\|_{L^2(\Om)}^2 + \|\theta^h\|_{L^\infty(\Om)}\right]~dt\right\}\,.\]
Applying Gr\"onwall's inequality together with \eqref{theta_int_bnd} and 
Lemma \ref{lemma:basich1} then gives
\begin{equation}\label{eq:a-final}
	\Ah(t) \leq  \tilde C\left[1 + \Bk(t) + \Dk(t)^\frac{1}{2}\right].
\end{equation}

\subsubsection*{Bound for $\Bk$}
To bound the $\Bk$ functional we define the test functions 
$$
	\psi^h (x,t):=  L'(\theta^h(x,t))\,,
$$
where $L$ is given in {\bf [A1]}.
Using $\psi^h$ as test function in \eqref{fem:tempeq}, with $\mu\equiv \bar \mu$, 
integrating in time, and rearranging yield
\begin{eqnarray*}
	&&\nqquad\nqquad   \int_\Om L(\theta^h(t)) ~dx
	+\frac{1}{h}\sum_{i=1}^{N-1} \int_0^t  G_i(\tau^h)\jump{L'(\theta^h)}_i^2~ds
	\nonumber \\
	&& \nqquad= \int_0^t \int_\Om  \Fh u^h_x L'(\theta^h)~dxds  + \int_\Om L(\theta^h)(0)~dx=: I_1 + I_2\,.
\end{eqnarray*}
Using the Cauchy-Schwarz inequality and the lower bounds in  \eqref{tau_bnd}-\eqref{theta_bnd}, gives
\begin{eqnarray*}
	I_1 &\leq& \tilde C\Big\{\int_0^t \int_\Om  |u^h_x|^2 |\theta^h|^{\beta+1}~dxds
	+   \int_0^t \int_\Om |\Fh|^2|\theta^h|^{\beta+1}~dxds\Big\}\nonumber  \\
	&\leq&\tilde C\Big\{\int_0^t \int_\Om  |u^h_x|^2 |\theta^h|^{\beta +\frac{3}{2}}~dxds 
	+ \int_0^t \int_\Om  |\theta^h|^2 |\theta^h|^{\beta +1}~dxds\Big\}\nonumber  \\
	&\leq&\tilde C\Big\{\Big[\sup_{s \in (0,t)}  \|\theta^h(s)\|_{L^\infty(\Omega)}^{2\beta +3})\Big]^\frac{1}{2}
	+\int_0^t \|\theta^h(s)\|_{L^\infty(\Omega)}\Big[\int_\Om  L(\theta^h)dx\Big] ~ds \Big\}\nonumber \\
	&\leq&\tilde C\Big\{1 + \Dk(t)^\frac{1}{2} + \int_0^t \|\theta^h(s)\|_{L^\infty(\Omega)}\Bk(s) ~ds\Big\}\,,
\end{eqnarray*}
where Lemma \ref{lem:sobolevestimate} is applied in the last inequality.
The $I_2$ term is bounded by the requirements on the initial data, whence 
\[\Bk(t) \leq \tilde C\Big\{1 + \Dk(t)^\frac{1}{2} + \int_0^t \|\theta^h(s)\|_{L^\infty}\Bk(s) ~ds\Big\}\,.\]
An application of Gr\"onwall's inequality, where we use the bound \eqref{theta_int_bnd}, yields 
\begin{equation}\label{eq:b-final}
	\Bk(t) \leq \tilde C\big[1+ \Dk(t)^\frac{1}{2}\big]\,.
\end{equation}

\subsubsection*{Bound for $\Dk$}

To bound $\Dk$ we define the test function $\psi^h\in Q_h(\Om)$ by
$$
\psi^h (x,t)= L'(\theta^h(x,t))_t = \kappa(\theta^h(x,t))\theta^h_t(x,t)\,.
$$
Using $\psi^h$ in the temperature scheme \eqref{fem:tempeq} and integrating in time give
\begin{equation*}
	\begin{split}
		&\int_0^t\int_\Om  \kappa(\theta^h)|\theta^h_t|^2 ~dxds 
		+\frac{1}{2h} \sum_{i =1}^{N-1} \int_0^tG_i(\tau^h)
		\big(\jump{L'(\theta^h)}_i^2\big)_t~ds \\
		&\qquad\qquad\qquad = \int_0^t \int_\Om  \Fh u_x^h \big(L'(\theta^h)\big)_t ~dxds\,.
	\end{split}
\end{equation*}
Integrating by parts in time and rearranging gives 
\begin{eqnarray}
	&&\nquad \frac{1}{2h}\sum_{i=1}^{N-1}  G(\tau^h) \jump{L'(\theta^h)}_i^2
	+\int_0^t \int_\Om  \kappa(\theta^h)|\theta^h_t|^2 ~dxds
	\nonumber \\
	&&\nquad =\int_0^t\! \int_\Om  \Fh u_x^h \big(L'(\theta^h)\big)_t dxds 
	+ \frac{1}{2h}\sum_{i=1}^{N-1} \Big[G_i(\tau^h)\jump{L'(\theta^h)}_i^2\Big](0)
	\nonumber \\
	&& 
	+\frac{1}{2h} \sum_{i=1}^{N-1}\int_0^t  \big[G_i(\tau^h)\big]_t\, \jump{L'(\theta^h)}_i^2~ds 
	 =: J_1 + J_2 + J_3\,.\label{eq:d-basic}
\end{eqnarray}
Using the Cauchy-Schwarz inequality with  parameter $\epsilon$ 
together with the pointwise bounds \eqref{tau_bnd}, \eqref{theta_bnd} on $\tau$ and $\theta$
(and the assumption $\beta>0$), we get
\[\begin{split}
	J_1 &\leq \epsilon \int_0^t \int_\Om \kappa(\theta^h)|\theta^h_t|^2~dxds 
	+ \frac{\tilde C}{\epsilon} \int_0^t \int_\Om \kappa(\theta^h)|\Fh|^2 |u_x^h|^2~dxds \\
	&\leq \epsilon \Dk(t) + \frac{\tilde C}{\epsilon}\int_0^t \|\kappa(\theta^h)\|_{L^\infty} 
	\left[\int_\Om |\theta^h|^4 +  |\Fh|^4 ~dx\right]~dt \\
	&\leq \epsilon \Dk(t) 
	+ \frac{\tilde C}{\epsilon}\Big[\sup_{s \in (0,t)}\int_\Om  L(\theta^h)~dx\Big]
	\left[\int_0^t  \|\theta^h\|_{L^\infty(\Omega)}^2~ds\right]\\
	&\quad + \frac{\tilde C}{\epsilon} \Big[\sup_{s \in (0,t)}\|
	\kappa(\theta^h)\|_{L^\infty(\Omega)}\Big]\!\!
	\Big[\sup_{s \in (0,t)}\int_\Om  |\Fh|^2~dx \Big]\!\!
	\left[\int_0^t \|\Fh\|_{L^\infty(\Omega)}^2~ds\right] \\
	&\leq \epsilon \Dk(t) + \frac{\tilde C}{\epsilon}\Bk(t)\int_0^t  \|\theta^h\|_{L^\infty(\Omega)}^2~ds \\
	&\quad + \frac{\tilde C}{\epsilon}\Big[\sup_{s \in (0,t)}
	\|\kappa(\theta^h)\|_{L^\infty(\Omega)}\Big]\!\!
	\left[\Ah(t) + \Bk(t) \right]\!\!\left[\int_0^t \|\Fh\|_{L^\infty(\Omega)}^2~ds\right].
\end{split}\]
Applying Lemma \ref{lemma:magicf} and the bounds \eqref{eq:a-final} and  \eqref{eq:b-final},
we obtain 
\[\begin{split}
	J_1 &\leq \epsilon \Dk(t)  + \frac{\tilde C}{\epsilon}\Bk(t)\left[1 + \Bk(t)^\frac{1}{2}\right] \\
	&\quad + \frac{\tilde C}{\epsilon}\Big[\sup_{s \in (0,t)}
	\|\kappa(\theta^h)\|_{L^\infty(\Omega)}\Big]\!\!
	\Big[\Ah(t) + \Bk(t)\Big]\!\!\Big[1+ \Ah(t)^\frac{1}{2} \Big]\\
	&\leq \epsilon \Dk(t) + \frac{\tilde C}{\epsilon}\left[1+ \Dk(t)^\frac{3}{4}\right]\!\!
	\Big[ 1+\sup_{s \in (0,t)}\|\kappa(\theta^h)\|_{L^\infty(\Omega)}\Big]\,.
\end{split}\]
To bound the last factor we apply Lemma \ref{lem:sobolevestimate}:
\[\sup_{s \in (0,t)}\|\kappa(\theta^h)\|_{L^\infty(\Omega)} 
=\tilde C \!\!\sup_{s \in (0,t)} 
\left[  \|\theta^h\|_{L^\infty(\Omega)}^{2\beta + 3}\right]^\frac{\beta}{2\beta + 3} 
\leq \tilde C\left[ 1 + \Dk(t)^\frac{\beta}{2\beta + 3}\right].\]
Thus, 
\begin{equation}\label{eq:d-i1}
	J_1 \leq \epsilon \Dk(t) + \frac{C}{\epsilon}\Dk(t)^\frac{10\beta+9}{4(2\beta+3)}\,,
\end{equation}
where we note that $\frac{10\beta+9}{4(2\beta+3)}<1$, since $\beta <\frac{3}{2}$. 
Next, consider $J_2$: by the requirements on the initial data
\begin{equation}\label{eq:d-i2}
	J_2 = {\frac{1}{2h}}\sum_{i=1}^{N-1} \Big[G_i(\tau^h_0)\jump{L'(\theta^h_0)}_i^2\Big] 
	    \leq {\frac{C}{h}}\Big[\sum_{i=1}^{N-1}\jump{\theta^h_0}_i^2\Big] \leq C\|\theta_{0,x}\|_{L^2(\Om)}^2.
\end{equation}
To bound $J_3$ we first observe that for each $1\leq i\leq N$,
\[ \big(G_i(\tau^h)\big)_t = -\textstyle\frac{1}{2}G_i(\tau^h)^2\left[u_x^h|_{E_i} + u_x^h|_{E_{i+1}}\right]\,.\]
Writing $u^h_x=\frac{\tau^h}{\bar \mu}(\Fh +p^h)$ and using the positivity of $\theta^h$ together 
with both  bounds in \eqref{tau_bnd}, we get
\[|G_i(\tau^h)_t|  \leq \tilde C G_i(\tau^h)\|\Fh\|_{L^\infty(\Omega)}\,,\qquad\mbox{where $\tilde C>0$,}\]
and thus:
\begin{eqnarray}\label{eq:d-i3}
&& \nqquad J_3=\frac{1}{2h} \sum_{i=1}^{N-1}\int_0^t  \big(G_i(\tau^h)\big)_t\, \jump{L'(\theta^h)}_i^2~ds
\nonumber\\
&\leq& \frac{\tilde C}{h}\int_0^t \left[\sum_{i=1}^{N-1} G_i(\tau^h) \jump{L'(\theta^h)}_i^2\right]\|\Fh\|_{L^\infty} ~ds 
\nonumber\\
&\leq& \tilde C  \Big[\sup_{s\in(0,t)}\sum_{i=1}^{N-1}  \frac{G_i(\tau^h)}{2h} \jump{L'(\theta^h)}_i^2\Big]^\frac{1}{2}
\int_0^t \|\Fh\|_{L^\infty} \Big[\sum_{i=1}^{N-1} \frac{G_i(\tau^h)}{h}\jump{L'(\theta^h)}_i^2\Big]^\frac{1}{2}ds
\nonumber\\
&\leq&\tilde C \Dk(t)^\frac{1}{2}\cdot \Big[\int_0^t \|\Fh\|_{L^\infty(\Omega)}^2~ds\Big]^\frac{1}{2}
\cdot \Big[\frac{1}{h}\sum_{i=1}^{N-1} \int_0^t  G_i(\tau^h)\jump{L'(\theta^h)}_i^2~ds\Big]^\frac{1}{2}
\nonumber\\
&\leq& \tilde C \Dh(t)^\frac{1}{2}\left[1+ \Ah(t)^\frac{1}{2} \right]^\frac{1}{2} \Bk(t)^\frac{1}{2}
\leq \tilde C \left[1+\Dk(t)^\frac{7}{8}\right]\,,
\end{eqnarray}
where we have used Lemma \ref{lemma:magicf} together with the previously derived
bounds on $\Ah$ and $\Bk$.
Substituting \eqref{eq:d-i1}, \eqref{eq:d-i2}, and \eqref{eq:d-i3} into \eqref{eq:d-basic}
we obtain 
\[\Dk(t) \leq \tilde C\left[1+\Dk(t)^{\delta}\right]\,,\qquad 
\mbox{with $\delta:=\max\left\{\frac{7}{8},\, \frac{10\beta+9}{8\beta+6}\right\}<1$,}\]
such that 
\begin{equation*}
	\Dk(t) \leq \tilde C\,.
\end{equation*}
Recalling \eqref{eq:a-final}, \eqref{eq:b-final} we conclude that there is a $C>0$,
depending on $T$, the systems parameters, and initial data, but independent of $h$, and such that
\begin{equation*}
	\Ah(t) + \Bk(t) + \Dk(t) \leq C\,,\qquad \mbox{for all $t\in[0,T]$.}
\end{equation*}
This concludes the proof of Lemma \ref{lemma:gradient}.  \qed


\subsubsection*{Technical lemmas used in the proof}
\begin{Lemma}\label{lemma:fbound}
Assume {\bf [A1]}-{\bf [A2]}. Then there is a $C$ independent of $h$ such that
\begin{eqnarray*}
	\nqquad\nqquad\left| \int_{0}^t \int_{\Omega} \Fh(x,s) \phi(x,s) ~dxds \right| &\leq& 
	\nonumber\\
	&&\nqquad  \nqquad \nqquad \nqquad \nqquad 
	C\left\{\epsilon \left[1 +  \Ah(t)\right]
	+ \frac{1}{\epsilon}\int_{0}^t \left[\int_{\Omega}|\phi(x,s)|~dx\right]^2 ds\right\}\,,
\end{eqnarray*}
for all $\epsilon>0$, $t\in(0,T)$ and $\phi \in L^2(0,T;L^1(\Om))$.
\end{Lemma}
\bProof
For each $t\in(0,T)$ define the test functions $v^h$ by 
\[v^h(x,t) =  \int_{0}^x [\Pi^Q_h \phi](y,t) - \left[\frac{1}{|\Om|}\int_\Om \phi(z,t) ~dz\right]~dy\,.\]
Using $v^h$ as test function in \eqref{fem:momentumeq}, rearranging and integrating in time, give 
\begin{equation*}
	\begin{split}
	&\left|\int_{0}^t \!\int_{\Omega} \Fh \phi~dxds\right| 
	= \left|\int_{0}^t \!\int_{\Omega}  u^h_{t}v^h dxds 
	- \frac{1}{|\Om|}\int_{0}^t   \left[\int_\Om \phi ~dx\right]\!\!\left[\int_\Om \Fh dx \right]ds\right| \\
	& \leq \int_0^t\! \left[\int_\Om |u^h_t|~dx\right]\!\!
	\left[\int_\Om |\phi|~dx\right]ds 
	+ C\int_0^t\!\left[\int_\Om |\phi|~dx\right]\!\!
	\left[\int_\Om \big(|u^h_x| + \theta^h\big)~ dx\right]ds \\
	&\qquad \leq \epsilon \Ah(t) 
	+ \frac{C}{\epsilon}\int_0^t \left[\int_\Om |\phi|~dx\right]^2ds
	+\epsilon\int_0^t \int_\Om |u^h_x|^2~dxds + C\epsilon t\,,
	\end{split}
\end{equation*}
where we have applied the Cauchy-Schwarz inequality with parameter $\epsilon$, 
together with the bound \eqref{eq:energyconv}. 
Applying the estimate \eqref{energy_0} concludes the proof.
\qed

\begin{Lemma}\label{lem:sobolevestimate}
Assume {\bf [A1]}-{\bf [A2]}. Then there is a $C$ independent of $h$ such that
\begin{equation*}
	\left[ \sup_{s \in (0,t)} \|\theta^h(s)\|_{L^\infty(\Om)}^{2\beta + 3}\right]^\alpha 
	\leq C\left[ 1 + \Dk(t)^\alpha\right],
\end{equation*}
for all $0<\alpha < 1$.
\end{Lemma}
\bProof
Fix any $s \in (0,t)$. Since $\theta^h=\theta^h(\cdot,s)$ is piecewise constant, we have
\begin{equation*}
\begin{split}
	\theta^h L'(\theta^h)^2\big|_{E_i}
	&= \theta^h\big|_{E_i}\Big[L'(\theta^h)^2\big|_{E_j} 
	\pm \sum_{k=j}^{i-1} \jump{L'(\theta^h)^2}_k\Big]\\
	&= \theta^h\big|_{E_i}L'(\theta^h)^2|_{E_j} 
	\pm \theta^h\big|_{E_i}\sum_{k=j}^{i-1} \jump{L'(\theta^h)}_k\big[L'(\theta^h)\big|_{E_k} 
	+L'(\theta^h)\big|_{E_{k+1}} \big].
\end{split}
\end{equation*}
Multiplying through by $h$, summing over $j= 1, \ldots, N$, taking the maximum 
over $i$, and applying the Cauchy--Schwartz inequality, we deduce
\begin{equation*}
	\begin{split}
	\|\theta^h(s)\|_{L^\infty(\Om)}^{2\beta + 3}
	&\leq \|\theta^h(s)\|_{L^\infty(\Om)}\int_\Om L'(\theta^h(s))^2~dx  \\
	& \qquad + 2\|\theta^h(s)\|_{L^\infty(\Om)}^{\frac{2\beta +3}{2}}
		\left[\int_\Om \theta^h(s)~dx \right]^\frac{1}{2}
		\left[\sum_{k=1}^{N-1}\frac{1}{h}\jump{L'(\theta^h(s))}_k^2\right]^\frac{1}{2} \\
	&\leq \|\theta^h(s)\|_{L^\infty(\Om)}^{2\beta +2} 
		+ \tilde C\|\theta^h(s)\|_{L^\infty(\Om)}^{\frac{2\beta +3}{2}}
		\left[\sum_{k=1}^{N-1}\frac{1}{h}\jump{L'(\theta^h(s))}_k^2\right]^\frac{1}{2} \\
	&\leq \epsilon \|\theta^h(s)\|_{L^\infty(\Om)}^{2\beta +3}  + 
	\tilde C(\epsilon)\left[1 + \sum_{k=1}^{N-1}\frac{1}{h}\jump{L'(\theta^h(s))}_k^2\right],
\end{split}
\end{equation*}
where we have used the energy bound \eqref{eq:energyconv}. 
We choose $\epsilon$ suitably small, absorb the first term on  the right-hand 
side into the left-hand side, and take the supremum over time. Recalling the definition 
of $\Dk(t)$, using the lower bound \eqref{tau_bnd}${}_1$, and taking the $\alpha$th 
power of both sides, then yields
\begin{equation*}
	\begin{split}
 	\Big[ \sup_{s \in (0,t)} \|\theta^h(s)\|_{L^\infty(\Om)}^{2\beta + 3}\Big]^\alpha
		&\leq \tilde C\Big[ 1 + \sup_{s \in (0,t)}\sum_{k=1}^{N-1}\frac{1}{h}
		\jump{L'(\theta^h)(s)}_k^2\Big]^\alpha \leq C\Big[ 1 + \Dk(t)^\alpha\Big]\,.
	\end{split}
\end{equation*}\qed

\begin{Lemma}\label{lemma:magicf}
Assume {\bf [A1]}-{\bf [A2]}. Then there is a $C$ independent of $h$ such that
\begin{equation}\label{eq:magicf}
	\int_0^t  \|\Fh(s)\|^2_{L^\infty(\Omega)}~ds 
\leq C\left[1+ \Ah(t)^\frac{1}{2} \right]\,.
\end{equation}
and
\begin{equation}\label{eq:magictheta}
\int_0^t  \|\theta^h(s)\|^2_{L^\infty(\Omega)}~ds 
\leq C\left[ 1+ \Bk(t)^\frac{1}{2}\right]\,,
\end{equation}
\end{Lemma}
\bProof
Consider the $i$th element $E_i$ of the mesh $\mathbb{E}_h$ and let $\Fh_i$ denote the 
constant value of $\Fh$ on this element. We have 
\[|\Fh_i|^2 = |\Fh_j|^2 \pm \sum_{k=j}^{i}\jump{|\Fh|^2}_k\,,\]
and we proceed to multiply by $\sigma^\frac{1}{2}h$, integrate in time, and sum over $j$. 
Making use of the lower bound on $\tau^h$, \eqref{eq:energyconv}, \eqref{theta_int_bnd} 
and Lemma \ref{lemma:basich1} to bound $\int\int|\Fh|^2~dxds$ we obtain
\begin{eqnarray}\label{eq:juleaften}
	&&\nqquad\int_0^t |\Fh_i|^2~ds \leq \tilde C\left\{\int_0^t\int_\Om |\Fh|^2~dxds
	+ \int_0^t\sum_{k=1}^{N-1}\jump{\Fh}_k(\Fh_k + \Fh_{k+1})~ds\right\} \nonumber\\
	&&\qquad\leq \tilde C\left\{1 +\left[\int_0^t \int_\Om |\Fh|^2 ~dxds\right]^\frac{1}{2}
	\left[\int_0^t\sum_{k=1}^{N-1} h \frac{\jump{\Fh}_k^2}{h^2}~ds \right]^\frac{1}{2} \right\} \nonumber\\
	&&\qquad\leq 
	\tilde C\left\{1 + \left[\int_0^t\sum_{k=1}^{N-1} h  \frac{\jump{\Fh}_k^2}{h^2}~ds \right]^\frac{1}{2}\right\}\,,
\end{eqnarray}
where we have used the Cauchy-Schwarz inequality.
To bound the remaining term, let $x_j$ be an arbitrary interior node and let $v^h \in V^0_h(\Om)$ be such that 
$v^h(x_k)= 1$ and $v^h(x_l) = 0$, $\forall l \neq k$. Using $v^h$ as a test function in 
\eqref{fem:momentumeq} we have
\[\jump{\Fh}_k=-\int_{E_{k+1}}\nquad \Fh_{k+1}v_x^h~dy + \int_{E_k}\Fh_kv_x^h~dy 
= \int_{E_{k}\cup E_{k+1}}\nqquad u_t^hv^h ~dx \,.\]
such that
\[\frac{|\jump{\Fh}_k| }{h} \leq  \big(\Pi^Q_h |u_t^h|\big)\big|_{E_k} +\big(\Pi^Q_h |u_t^h|\big)\big|_{E_{k+1}}\,.\]
Since $k$ is arbitrary it follows that 
\[\sum_{k=1}^{N-1}h  \frac{\jump{\Fh}_k^2}{h^2}
\leq \tilde C\int_\Om \big(\Pi^Q_h|u^h_t|\big)^2\, dx\leq \tilde C\int_\Om|u^h_t|^2\, dx\,.\]
Using this in \eqref{eq:juleaften} yields
\[\int_0^t |\Fh_i|^2~ds \leq \tilde C\left[ 1 + \Ah(t)^\frac{1}{2}\right]\,,\]
and \eqref{eq:magicf} follows.
To establish \eqref{eq:magictheta} we argue similarly to get
\[|\theta_i^h|^2 \leq \left\{ \int_\Om |\theta^h|^2\, dx + \sum_{k=1}^{N-1}
\big(\theta^h_k + \theta^h_{k+1}\big)|\jump{\theta^h}_k|  \right\}\,.\]
Using \eqref{eq:energyconv} and \eqref{theta_int_bnd} to bound $\int\int|\theta^h|^2~dxds$,
together with $\eqref{tau_bnd}_1$, we obtain
\begin{eqnarray*}
	&&\int_0^t |\theta^h_i|^2~ds \leq \tilde C\left\{1
	+ \int_0^t\sum_{k=1}^{N-1}(\theta^h_k + \theta^h_{k+1})|\jump{\theta^h}_k|~ds\right\} \\
	&&\qquad\leq  \tilde C\left\{1 +\left[\int_0^t \int_\Om |\theta^h|^2 ~dxds\right]^\frac{1}{2}
	\left[\int_0^t\frac{1}{h}\sum_{k=1}^{N-1} \jump{\theta^h}_k^2~ds \right]^\frac{1}{2} \right\}\\
	&&\qquad\leq 
	\tilde C\left\{1 + \left[\frac{1}{h}\sum_{k=1}^{N-1} \int_0^t
	G_k(\tau^h)\jump{L'(\theta^h)}_k^2~ds \right]^\frac{1}{2}\right\}
	\leq \tilde C\left[ 1+ \Bk(t)^\frac{1}{2}\right] \,,
\end{eqnarray*}
where we have used that 
$\jump{L'(\theta^h)}_k=\kappa(\theta^*)\jump{\theta^h}_k\geq \tilde C\jump{\theta^h}_k$,
(for a $\theta^*$ between $\theta_k$ and $\theta_{k+1}$) due to the lower bound \eqref{theta_bnd}. 
\qed

\section{Proof of Lemma \ref{lemma:pointwise}}\label{sec:pointwise}
Throughout this section the assumptions {\bf [A1]} and {\bf [A2]} in Section \ref{sec:proof} are in force,
and we fix an arbitrary time $T>0$.
To simplify the notation we set
\[\omega^h (x,t):= \int_0^t \Fh(x,s)~ds + \Pi_h^Q \left[\int_0^x u^h_0(y)~dy\right]\,.\]
In what follows $\tilde C$, $C_1,\dots$ denote positive numbers that  depend on $T$, $|\Omega|$, 
the initial data, and the bounds \eqref{data_bnds}, but that are independent of $h$.
%
%
\begin{Lemma}\label{pws_lemma1}
There is a number $0<C_1<\infty$ such that
\[ |\omega^h(x,t)| \leq C_1\,,\qquad\mbox{for all $(x,t)\in\Omega\times [0,T]$.}\]
\end{Lemma}
\bProof
For fixed $t>0$ and $v^h\in V_h^0(\Omega)$, integrating \eqref{fem:momentumeq} in time gives
\begin{equation}\label{eq:omega-basic}
	\int_\Om \omega^h v^h_x~dx 
	= -\int_\Om u^hv^h ~dx\,.
\end{equation}
Using $v^h(x)=u^h(x,t)$ then yields
\[\int_\Om \omega^h u^h_x~dx = -\int_\Om |u^h|^2 ~dx\,.\]
Using this together with \eqref{fem:contequation} with $\phi^h(x) = \omega^h(x,t)$, gives\footnote{At this 
point the assumption of constant viscosity is used in an essential manner.}
\[\frac{d}{dt}\int_\Om \omega^h \tau^h~dx 
= \int_\Om \Fh \tau^h~dx + \int_\Om \omega^h u^h_x~dx 
= -\int_\Om K\theta^h + |u^h|^2~dx\,.\]
According to \eqref{eq:energyconv} we thus have, for each $t\geq 0$, that
\begin{equation}\label{uppr_lwr}
	\tilde C-t(K+2)\mathcal{E}(\theta^h_0, u^h_0)
	\leq \left(\int_\Om \tau^h \omega^h dx\right)\!\!(t)
	\leq \int_\Om \omega^h(x,0) \tau^h_0(x)~dx =: \tilde C.
\end{equation}
Next, by conservation of mass \eqref{eq:massconv} there are maps 
$t \mapsto z^\pm(t) \in \Om$ such that
\begin{equation}\label{eq:omega2}
	\omega^h(z^-(t),t) \leq \left(\int_\Om \tau^h \omega^h dx\right)(t) \leq \omega^h(z^+(t),t), \quad \forall t.
\end{equation}
Next, fix $t$, $h>0$, and $y \in \Om$ arbitrary. Let $E^- \in \mathbb{E}_h$ be the element
containing $z^-(t)$, let $E^y$ be the element containing $y$, and define
the interval $S^y \subseteq \Om$ by
$$
S^y := \left[\min\{z^-(t), y\}, \max\{z^-(t), y\}\right]\cup E^- \cup E^y\,.
$$
Then fix a $v^h \in V^0_h(\Om)$ by requiring 
that for each $x_i$, $i \in \Set{0,1,\ldots, N}$,
\begin{equation*}
	v^h(x_i) = 
	\begin{cases}
		1& \textrm{ if } x_i \in S^y, \\
		0& \textrm{ otherwise }.
	\end{cases}
\end{equation*}
Since $\omega^h$ is piecewise constant we have
\begin{equation*}
	\begin{split}
	&\omega^h(y,t) - \omega^h(z^-(t),t) = \pm\sum_{\{i: x_i \in S^y\}}\jump{\omega^h}_i \\
	&\qquad = \pm \sum_{i= 1}^{N-1}\jump{\omega^h}_iv^h(x_i) = \mp\int_\Om \omega^h v_x^h~dx 
	= \pm\int_\Om u^h v^h~dx,
\end{split}
\end{equation*}
where the last equality follows by \eqref{eq:omega-basic}. Rearranging, applying H\"older and the 
energy estimate \eqref{eq:energyconv}, using \eqref{eq:omega2} together with \eqref{uppr_lwr},
and finally taking the supremum in $y$, give
\[\sup_{y \in \Om} \omega^h(y,t) \leq \tilde C\,.\]
A similar argument using $z^+(t)$ yields $\inf_{y \in \Om} \omega^h(y,t) \geq -\tilde C$.\qed

%
%
\begin{Lemma}\label{pws_lemma2}
	There exists $C_2>0$ such that
	\[C_2^{-1} \leq \|\tau^h(t)\|_{L^\infty(\Omega)} \leq C_2\left[1 +\int_0^t\|\theta^h(s)\|_{L^\infty(\Omega)}~ds\right]
	\qquad\mbox{for $t\in [0,T]$}\,.\]
\end{Lemma}
\bProof
Since $\tau^h_t(\cdot,t),\, u^h_x(\cdot,t)\in Q_h(\Om)$ at each time, \eqref{fem:contequation}
gives $\tau^h_t = u^h_x$. Thus, with $A:=(\tau^h)^{-1}\exp(\bar\mu^{-1}\omega^h)$ and $B:=A^{-1}$,
we get
\[A_t=-\frac{p(\theta^h, \tau^h)}{\bar\mu}\, A\leq 0\,,\qquad B_t=\frac{p(\theta^h, \tau^h)}{\bar\mu}\, B
= \frac{K}{\bar\mu }\theta^h \exp(-\frac{\omega^h}{\bar\mu})\,.\]
Integration in time, together with \eqref{data_bnds} and Lemma \ref{pws_lemma1}, give
\[A(x,t)\leq A_0(x)\leq \tilde C\quad\Rightarrow\quad \tau^h(x,t)\geq C_2^{-1}\,,\]
and
\[B(x,t)\leq  B_0(x)+\tilde C\int_0^t\!\! \|\theta^h(s)\|_{\infty}ds
\,\,\Rightarrow\,\, \tau^h(x,t)\leq C_2\left[1 +\int_0^t\!\!\|\theta^h(s)\|_{\infty}~ds\right].\]
\qed

%
%
\begin{Lemma}\label{pws_lemma3}
	There exists $C_3>0$ such that
	\begin{equation}\label{lower_theta}
		\left\|\frac{1}{\theta^h(t)}\right\|_{L^\infty(\Omega)} \leq C_3   \qquad\mbox{for $t\in [0,T]$}\,.
	\end{equation}
\end{Lemma}
\bProof
Apply \eqref{fem:tempeq} with $\psi^h = -M(\theta^h)^{-M-1}$, where $M>0$.
A calculation (completing the square in $u^h_x$) shows that 
\[\frac{d}{dt}\int_\Om (\theta^h)^{-M}~ dx - \frac{M}{h}\sum_{i=1}^{N-1}
	G_i(\tau^h)\jump{L'(\theta^h)}_i\!\!\left\llbracket \frac{1}{(\theta^h)^{M+1}}\right\rrbracket_i  
	\leq \int_\Om \frac{MK^2dx}{4\bar \mu\tau^h (\theta^h)^{M-1}}\,.\]
Since $L'=\kappa>0$ and $M>0$ the sum on the left is negative, and H\"older gives
\begin{equation*}
	\frac{d}{dt}\left[ \left\|\frac{1}{\theta^h(t)}\right\|_M^M \right]
	\leq \frac{MK^2}{4\bar \mu}\int_\Om \frac{dx}{\tau^h (\theta^h)^{M-1}} 
	\leq \tilde CM\left\|\frac{1}{\tau^h(t)}\right\|_M\left\|\frac{1}{\theta^h(t)}\right\|_M^{M-1}\,,
\end{equation*}
where $\|\cdot\|_M$ denotes the ${L^M(\Om)}$-norm. Applying 
Lemma \ref{pws_lemma2} we obtain
\begin{equation*}
	\frac{d}{dt}\left[\left\|\frac{1}{\theta^h(t)}\right\|_M \right]
	\leq \tilde C\left\|\frac{1}{\tau^h(t)}\right\|_M \leq \tilde C\,.
\end{equation*}
Integrating in time, applying $\eqref{data_bnds}_2$, and sending $M\uparrow\infty$, yield \eqref{lower_theta}.\qed

%
%
\begin{Lemma}\label{pws_lemma4}
There exists $C_4>0$ such that
\begin{equation}\label{eq:middle-theta}
	\int_0^T \|\theta^h(s)\|_\infty~ds \leq C_4\,.
\end{equation}
\end{Lemma}
\bProof
Fix $x,\, y \in \Om$. As $\theta^h(\cdot,t)\in Q_h(\Omega)$,  
Lemma \ref{pws_lemma3} yields (suppressing $t$)
\[\theta^h(x)^\frac{1}{2} 
	= \theta^h(y)\theta^h(x)^{-\frac{1}{2}} + \sum_{i=j}^{k-1}\frac{\jump{\theta^h}_i}{\theta^h(x)^{\frac{1}{2}}}
	\leq \tilde C\theta^h(y) + \sum_{i=1}^{N-1}\frac{|\jump{\theta^h}_i|}{\theta^h(x)^{\frac{1}{2}}}\,,\]
where $x \in E_k$ and $y \in E_j$. Integrating in $y$, using \eqref{eq:energyconv}, and choosing $x$ 
such that $\|\theta^h(t)\|_{L^\infty(\Omega)}=\theta^h(x,t)$, yield
\[\|\theta^h(t)\|_{L^\infty(\Omega)}^\frac{1}{2}\leq 
\tilde C\left[1+\sum_{i=1}^{N-1}\frac{|\jump{\theta^h}_i|}{\|\theta^h(t)\|_{L^\infty(\Omega)}^\frac{1}{2}}\right]
\leq \tilde C\left[1+\sum_{i=1}^{N-1}\frac{|\jump{\theta^h}_i|}{(\theta^h_{i*})^\frac{1}{2}}\right]\,,\]
where $\theta^h_{i*}$ is as in \eqref{entrp}.
Multiplying and dividing by suitable terms in the latter sum, using $\kappa(\theta)\propto \theta^\beta$,
applying Cauchy-Schwarz, and squaring, give
\begin{equation}\label{intermed}
	\|\theta^h(t)\|_{L^\infty(\Omega)} \leq
	\tilde C\left\{1+\left[\sum_{i=1}^{N-1}\frac{G_i(\tau^h)
	\kappa(\theta^h_{i\dagger})}{h(\theta^h_{i*})^2}\jump{\theta^h}_i^2\right]\!\!
	\left[\sum_{i=1}^{N-1} \frac{h\theta^h_{i*}}{G_i(\tau^h)(\theta^h_{i\dagger})^\beta}\right]\right\},
\end{equation}
where $\theta^h_{i\dagger}$ is as in \eqref{entrp}.
Let the two inner sums on the right-hand side of \eqref{intermed} be $\mathfrak A(t)$ 
and $\mathfrak B(t)$, respectively.
To bound $\mathfrak B(t)$ we recall that $G_i(\tau^h)= 2\big(\tau^h|_{E_i} + \tau^h|_{E_{i+1}}\big)^{-1}$, and use that 
$\theta^h_{i*}\leq\theta^h|_{E_i}+\theta^h|_{E_{i+1}}$, while $\theta^h_{i\dagger}$ is between 
$\theta^h|_{E_i}$ and $\theta^h|_{E_{i+1}}$. Lemma \ref{pws_lemma3}, \eqref{eq:energyconv}, 
and Lemma \ref{pws_lemma2}, then give
\[\mathfrak B(t)\leq \tilde C\|\tau^h(t)\|_{L^\infty(\Omega)} \leq \tilde C\Big[1 +\int_0^t\|\theta^h(s)\|_{L^\infty(\Omega)}~ds\Big]\,.\]
Squaring both sides in \eqref{intermed} thus gives
\[\|\theta^h(t)\|_{L^\infty(\Omega)} \leq \tilde C\left[ 1+ \mathfrak A(t)\cdot\Big[1 +\int_0^t\|\theta^h(s)\|_{L^\infty(\Omega)}~ds\Big]\right]\,.\]
Applying the Gr\"onwall inequality and recalling, by \eqref{eq:entropy}, that $\int_0^t \mathfrak A(s)\, ds \leq C_0$,
we obtain \eqref{eq:middle-theta}.\qed

\noindent
Combing Lemma \ref{pws_lemma1} - Lemma \ref{pws_lemma4} completes the proof 
of Lemma \ref{lemma:pointwise}.


\end{document}